\documentclass[preprint,10pt]{elsarticle}

\usepackage{amssymb}
\usepackage{amsthm}
\usepackage{graphicx}
\usepackage{epstopdf}
\usepackage{amsmath}
\usepackage{amsfonts}
\usepackage{natbib}
\usepackage[all]{xy}
\usepackage{amscd}
\usepackage{color}

\makeatletter

\newcommand{\Rmnum}[1]{\expandafter\@slowromancap\romannumeral #1@}
\makeatother

\newtheorem{theorem}{Theorem}[section]
\newtheorem{proposition}{Proposition}[section]
\newtheorem{corollary}{Corollary}[section]
\newtheorem{definition}{Definition}[section]
\newtheorem{example}{Example}[section]
\newtheorem{lemma}{Lemma}[section]
\newtheorem{remark}{Remark}[section]
\numberwithin{equation}{section}

\journal{}

\begin{document}

\begin{frontmatter}

\title{Correlation entropy of free semigroup actions}

\author[add1]{Xiaojiang Ye\corref{}}
\ead{yexiaojiang12@163.com}
\author[add1]{Yanjie Tang\corref{}}
\ead{yjtang1994@gmail.com}
\author[add1]{Dongkui Ma\corref{cor1}}
\cortext[cor1]{Corresponding author}
\ead{dkma@scut.edu.cn}
\address[add1]{School of Mathematics, South China University of Technology,\\
Guangzhou, 510640, China}

\begin{abstract}
This paper introduces the concepts of correlation entropy and local correlation entropy for free semigroup actions on compact metric space, and explores their fundamental properties. Thereafter, we generalize some classical results on correlation entropy and local correlation entropy to apply to free semigroup actions. Finally, we establish the relationship between topological entropy, measure-theoretic entropy, correlation entropy, and local correlation entropy for free semigroup actions under various conditions.
\end{abstract}


\end{frontmatter}

\section{Introduction}
In a classical dynamical system $(X,f)$, where $X$ is a compact metric space with metric $d$ and $f: X \rightarrow X$ is a continuous transformation,
a point $x$ is recurrent if, for any neighborhood $U$ of $x$, there exist infinitely many indices $n$ such that $f^{n}(x) \in U$. The topological version of the famous Poincar$\acute{e}$ recurrence theorem states that almost every point is recurrent for every $f$-invariant Borel finite measure. Due to the continuity of $f$, for any given $\varepsilon$ and any recurrent point $x$, there exist infinitely many pairs of indices $i \neq j$ such that $d(f^i(x), f^j(x))< \varepsilon$. These pairs are referred to as recurrences. In \cite{EKR}, Eckmann, Kamphorst and Ruelle explored recurrences by recurrence plots, a white-and-black square image with black pixels representing recurrences. Further insights into recurrence quantification analysis are available in literature such as \cite{MRTK,WM,ZW}. Building upon the investigation of recurrence, the correlation sum $C(x,\varepsilon,d,n)$ is defined as
$$
C(x,\varepsilon,d,n):=\frac{1}{n^2} \sharp \left\{(i,j):0 \leq i,j \leq n-1, d(f^i(x),f^j(x))\leq \varepsilon\right\},
$$
where $\varepsilon$ represents the threshold distance. Literature concerning the correlation sum can be found in references such as \cite{ABDGHW,DK,GMS}. \textcolor{red}{A fundamental result regarding the correlation sum states that in an ergodic dynamical system $(X,f,\mu)$, there exists a countable subset $Q \subset \mathbb{R}$ such that for any $\varepsilon \notin Q$, the following convergence holds for almost everywhere $x \in X$:}
\begin{equation}\label{1.1}
\lim_{n \rightarrow +\infty}C(x,\varepsilon,d,n)=\int_X \mu (B_d(x,\varepsilon))d\mu(x),
\end{equation}
where $B_d(x,\varepsilon)$ is the $\varepsilon$-neighborhood of $x$, and $\int_X \mu (B_d(x,\varepsilon))d\mu(x)$ is referred to as the correlation integral. A detailed proof of this statement can be found in \cite{ABDGHW,MS,P,PT,S1996}.\par
In the pursuit of numerical estimation method for generalized entropies, Takens introduced $q$-correlation entropy, expanding upon the concepts of correlation sum and correlation integral in \cite{T}. If $q \neq 1$, the $q$-correlation entropy is defined as
$$
\begin{aligned}
\overline{h}_{cor}(f,\mu,q)&:=\lim_{\varepsilon \rightarrow 0}\limsup_{k \rightarrow +\infty}-\frac{1}{(q-1)k}\log
\int_X\mu(B_{d_k}(x,\varepsilon))^{q-1}d\mu(x), \\
\underline{h}_{cor}(f,\mu,q)&:=\lim_{\varepsilon \rightarrow 0}\liminf_{k \rightarrow +\infty}-\frac{1}{(q-1)k}\log
\int_X\mu(B_{d_k}(x,\varepsilon))^{q-1}d\mu(x),
\end{aligned}
$$
if $q=1$, the $1$-correlation entropy is defined as
$$
\begin{aligned}
\overline{h}_{cor}(f,\mu,1)&:=\lim_{\varepsilon \rightarrow 0}\limsup_{k \rightarrow +\infty}-\frac{1}{k}
\int_X\log\mu(B_{d_k}(x,\varepsilon))d\mu(x), \\
\underline{h}_{cor}(f,\mu,1)&:=\lim_{\varepsilon \rightarrow 0}\liminf_{k \rightarrow +\infty}-\frac{1}{k}
\int_X\log\mu(B_{d_k}(x,\varepsilon))d\mu(x),
\end{aligned}
$$
where $d_k$ represents the Bowen metric and $B_{d_k}(x,\varepsilon)$ denotes the Bowen dynamical ball. Following \cite{T}, $\check{S}$pitalsk$\acute{y}$ \cite{S2018} defined local correlation entropy as follows
$$
  \begin{aligned}
   \overline{h}_{cor}(f,x):=\lim_{\varepsilon \rightarrow 0} \limsup_{k \rightarrow +\infty} - \frac{1}{k} \log \liminf_{n \rightarrow
   +\infty} C(x,\varepsilon,d_k,n), \\
   \underline{h}_{cor}(f,x):=\lim_{\varepsilon \rightarrow 0} \liminf_{k \rightarrow +\infty} - \frac{1}{k} \log \limsup_{n
   \rightarrow +\infty} C(x,\varepsilon,d_k,n).
  \end{aligned}
$$
In \cite{V}, Verbitskiy demonstrated that $1$-order correlation entropy equals measure-theoretic entropy. Moreover, if $\mu$ is an invariant $f$-homogeneous Borel probability measure, meaning that Borel probability measure $\mu$ is invariant and satisfies that for any
$\varepsilon>0$, there exists $\delta >0$ and $c>0$ such that $\mu(B_{d_k}(y,\delta)) \leq c\mu(B_{d_k}(x,\varepsilon))$ for all $x$,
$y\in X$ and $k\in \mathbb{N}$, then $q$-order correlation entropy, measure-theoretic entropy and topological entropy coincide
for any $q\in \mathbb{R}$. Based on formula (\ref{1.1}), we have
$$
\overline{h}_{cor}(f,x)=\overline{h}_{cor}(f,\mu,2), \qquad \underline{h}_{cor}(f,x)=\underline{h}_{cor}(f,\mu,2), \quad \mu-a.e. x.
$$
Thus, a connection exists among different entropies if $\mu$ is an invariant $f$-homogeneous Borel probability measure. This connection is valuable as it obviates the need for partition consideration when computing measure-theoretic entropy. \par

Recently, there has been growing interest in studying free semigroup actions, which allow systems to adapt dynamically over time to accommodate inevitable experimental errors. Ghys et al. \cite{GLW} proposed a definition for topological entropy applicable to finitely generated pseudo-groups of continuous maps, sparking significant interest in free semigroup actions on compact metric spaces. For instance, Bufetov \cite{B1999} introduced the concept of topological entropy for free semigroup actions, whereas Bi$\acute{s}$ \cite{B2004} introduced the entropies for a semigroup of maps using alternative methods. Partial variational principles linking measure-theoretic entropy and topological entropy for free semigroup actions were investigated in \cite{CRV1} and \cite{LMW}. Carvalho et al. proposed a novel definition for the measure-theoretic entropy of free semigroup actions in \cite{CRV2} to establish these principles. In contrast to studies involving the entire space of free semigroup actions, Ju et al. \cite{JMW} explored the topological entropy, while Xiao and Ma \cite{XM,XM2} investigated the topological pressure of free semigroup actions on non-compact sets. Given the interconnectedness among topological entropy, measure-theoretic entropy, correlation entropy, and local correlation entropy in classical dynamical systems, it is natural to inquire whether such relationships persist in free semigroup actions, although limit literature exists on this topic. Therefore, this paper aims to investigate this relationship. \par

Throughout this paper, we focus on a compact metric space $(X,d)$ and a Borel probability measure $\mu$ with full support. Let $G$ be the free semigroup with $m$ generators $\{f_1,f_2,\cdots,f_m\}$ acting on $X$, where each $f_i$ denotes a continuous self-map on $X$ with $i \in \{1,2,\cdots,m\}$. We use $h_{\mu}(G)$ to denote the measure-theoretic entropy, $h_{top}(G)$ to denote the topological entropy and $h(\omega,x)\left(H(\omega,x)\right)$ to denote the lower(upper) local entropy of free semigroup actions, respectively. \par
We introduce the concepts of correlation sum, upper(lower) local correlation entropy and $q$-order upper(lower) correlation entropy for free semigroup actions, denoted by $C(G,x,\varepsilon,\omega,k,n)$, $\overline{h}_{cor}(G,x)$($\underline{h}_{cor}(G,x)$), $\overline{h}_{cor}(G,\mu,q)$($\underline{h}_{cor}(G,\mu,q)$) respectively (details can be found in Section 3). Subsequently, we establish the relationship among different entropies. Specifically, Theorem \ref{th1.1} corresponds to the case of $q=2$, Theorem \ref{th1.2} to $q=0$, Theorem \ref{th1.3} to $q \geq 1$, while Theorem \ref{th1.4} and Theorem \ref{th1.5} address the cases of $0 \leq q \leq 1$ and $q \leq 0$, respectively. \par
We now proceed to present our findings.

\begin{theorem}\label{th1.1}
Let $G$ be a free semigroup acting on a compact metric space $X$, $\mu$ be $G$-ergodic Borel probability measure on $X$, \textcolor{red}{then there exists a countable subset $Q \subset \mathbb{R}$ such that for any $\varepsilon \notin Q$, $\omega \in \Sigma_m^+$ and $k \in \mathbb{N}$,}
\begin{equation}\label{eq1.2}
\lim_{n \rightarrow +\infty} C(G,x,\varepsilon,\omega,k,n) =\int_X \mu (B_{\omega,k}^G(x,\varepsilon))d\mu(x), \qquad \mu-a.e. \quad x \in X,
\end{equation}
where $B_{\omega,k}^G(x,\varepsilon)$ is generalized Bowen dynamical ball (see details in Section 2).
Furthermore,
$$
\overline{h}_{cor}(G,\mu,2)=\overline{h}_{cor}(G,x), \quad
\underline{h}_{cor}(G,\mu,2)=\underline{h}_{cor}(G,x), \qquad \mu-a.e. \quad x \in X.
$$
\end{theorem}
\begin{theorem}\label{th1.2}
Let $G$ be a free semigroup acting on a compact metric space $X$, $\mu$ Borel probability measure on $X$. If $\mu$ satisfies the weak entropy-doubling condition of free semigroup actions, then
$$
h_{top}(G)=h_{cor}(G,\mu,0).
$$
\end{theorem}

\begin{remark}
Theorem \ref{th1.1} and \ref{th1.2} are generalizations of classical results.
\end{remark}

\begin{theorem}\label{th1.3}
Let $G$ be a free semigroup acting on a compact metric space $X$, $\mu$ $G$-ergodic Borel probability measure on $X$. If $\mu$ satisfies $h_{\mu}(G) < +\infty$ and the limit process of $h(\omega,x)$ is uniformly about $x$ for almost everywhere $\omega$, then for any $q \geq 1$,
$$h_{\mu}(G)=h_{cor}(G,\mu,q).$$
In particularly,
$$
h_{\mu}(G)=h_{cor}(G,\mu,2)=h_{cor}(G,x), \qquad \mu-a.e. \quad x \in X.
$$
\end{theorem}

\begin{remark}
	In Theorem \ref{th1.3}, the condition regarding $h(\omega,x)$ is necessary, but its optimality is unknown. In general, without the condition of $h(\omega,x)$, the Theorem \ref{th1.3} does not hold. Examples demonstrating this can be found in \cite{S2018}.
\end{remark}

\begin{theorem}\label{th1.4}
Let $G$ be a free semigroup acting on a compact metric space $X$, $\mu$ be Borel probability measure on $X$. If $h(\omega,x) \geq h_{top}(G)$ almost everywhere, then for any $0 \leq q \leq 1$,
$$
h_{top}(G)=h_{cor}(G,\mu,q).
$$
Moreover, if $\mu$ is $G$-invariant, then
$$
h_{top}(G)=h_{cor}(G,\mu,1)=h_{\mu}(G).
$$
\end{theorem}
\begin{theorem}\label{th1.5}
Let $G$ be a free semigroup acting on a compact metric space $X$, $\mu$ be $G$-invariant Borel probability measure on $X$ satisfying the weak entropy-doubling condition of free semigroup actions. If the limit process of $H(\omega,x)$ is uniformly about $(\omega,x)$ and $H(\omega,x) \leq h_{\mu}(G)$ almost everywhere, then for any $q \leq 0$,
$$
h_{cor}(G,\mu,q) = h_{\mu}(G).
$$
In particularly,
$$
h_{top}(G)=h_{cor}(G,\mu,0)= h_{\mu}(G).
$$
\end{theorem}

\begin{remark}
Theorem \ref{th1.3}, \ref{th1.4} and \ref{th1.5} were first proposed even in classical dynamical systems and were inspired by the examples in chapter 2 of \cite{V}.
\end{remark}

The paper is structured as follows: Section 2 presents the necessary preliminaries, while Section 3 introduces the concepts of correlation entropy and local correlation entropy of free semigroup actions, along with an examination of their properties. The proofs of Theorems \ref{th1.1} and \ref{th1.2} are provided in Sections 4 and 5, respectively. Section 6 is dedicated to the proofs of Theorems \ref{th1.3}, \ref{th1.4}, and \ref{th1.5}.

\section{Preliminaries}

\subsection{The free semigroup actions on compact metric space}
Let $(X,d)$ be a compact metric space, and let $G$ be a free semigroup generated by $G_{\ast}:=\{f_1, f_2, \cdots, f_m\}$, where each $f_i$ is a
continuous self-map on $X$ for all $1\leq i\leq m$. Given a vector $p=(p_1, p_2, \cdots, p_m)$ with $\sum_{i=1}^{m}p_i=1$ and $p_i>0$
for all $1\leq i \leq m$, there exists a symbol space $\Sigma_m^+:= \{1, 2, \cdots, m\}^\mathbb{N}$ with a Bernoulli probability
measure $\mathbb{P}$ generated by the vector $p$. Let $\sigma: \Sigma_m^+ \longrightarrow \Sigma_m^+$ be the shift operator defined by
$\sigma (i_1, i_2, \cdots)=(i_2, i_3, \cdots)$. For any $\omega\in \Sigma_m^+$, denoted as $\omega=(i_1, i_2, \cdots)$, we define
$$
  f_{\omega,n}(x):= \begin{cases}
                    x, & n=0 \\
                    f_{i_n}\circ f_{i_{n-1}}\circ \cdots \circ f_{i_1}(x), & n\geq 1.
                   \end{cases}
$$
Thus, the orbit of $x$ under the free semigroup actions is defined as
$$
Orb(x,G):=\{f_{\omega,n}(x): \forall \omega \in \Sigma_m^+, n \geq 0\}.
$$
We refer to $(X,G)$ as the free semigroup action. Given $\omega=(i_1, i_2, \cdots)$ and $k \geq 1$, one could define the generalized Bowen metric as
$$d_{\omega, k} ^G (x, y):=\max \{d(f_{\omega, i}(x), f_{\omega, i}(y)) : 0 \leq i \leq k-1 \}
$$
and the generalized Bowen dynamical ball as
$$
B_{\omega,k}^G(x,\varepsilon):=\{y: d_{\omega,k}^G(x,y) \leq \varepsilon\}.
$$
Let $F: \Sigma_m^+ \times X \longrightarrow \Sigma_m^+ \times X$ be the skew product transformation defined as follows
$$
  F(\omega, x)=(\sigma (\omega), f_{\omega,1}(x)).
$$\par

Here, we revisit certain terminologies originating from random dynamical systems (cf. \cite{K,KL} for detailed exposition) and introduce specific constraints to adapt them for free semigroup actions. A Borel probability measure $\mu$ is called $G$-invariant if $\mathbb{P} \times \mu$ is invariant with respect to $F$. Similarly, A Borel probability measure $\mu$ is called $G$-ergodic if $\mathbb{P} \times \mu$ is both invariant and ergodic with respect to $F$. Subsequently, we present a necessary theorem concerning ergodicity in the random dynamical systems. While initially established by \cite{K1951}, the theorem underwent generalization by \cite{O}, and Kifer \cite{K} provided an alternative proof methodology. We apply this theorem to the context of free semigroup actions for our convenience.

\begin{theorem}\label{th2.1}
Let $G$ be a free semigroup acting on a compact metric space $X$, $\mu$ be a Borel probability measure on $X$. Then $\mathbb{P}\times\mu$ is ergodic with respect to skew product transformation if and only if for any measurable subset $A \subseteq X$, if $\sum_{i=1}^m p_{i} \chi_{_{A}}(f_i(x))=\chi_{_{A}}(x)$, then $\mu(A)=0$ or $1$.
\end{theorem}  \par

Given any integer $t \geq 1$, $G^t$ can be defined as a free semigroup generated by
$$
G_{\ast}^t:=\{g_1, g_2, \cdots, g_{m^t}:g_i=f_{i_t}\circ f_{i_{t-1}}\circ \cdots \circ f_{i_1}, 1 \leq i\leq m^t, f_{i_s}\in G_{\ast},
1\leq s\leq t\}.
$$
Additionally, a one-to-one transformation $\tau$ is defined as follows
$$
\begin{aligned}
\tau : \{1,2,\cdots, m^t\} &\longrightarrow \{1,2,\cdots,m\}^t\\
j&\mapsto(i_1,i_2,\cdots,i_t),
\end{aligned}
$$
and for convenience, the following transformation is still denoted as $\tau$
$$
\begin{aligned}
\tau : \Sigma_{m^t}^+ &\longrightarrow \Sigma_m^+\\
(j_1,j_2,\cdots)&\mapsto (\tau(j_1),\tau(j_2),\cdots).
\end{aligned}
$$
Thus, a probability vector $p^t$ on $\{1,2,\cdots,m^t\}$ corresponding to $p$ can be defined as $p^t(j)=p(i_1)p(i_2)\cdots p(i_t)$ where
$1\leq j\leq m^t$, $\tau(j)=(i_1,i_2,\cdots,i_t)$. Additionally, there exists a symbol space $\Sigma_{m^t}^+:=\{1,2,\cdots,m^t\}^{\mathbb{N}}$ with a Bernoulli
probability measure $\mathbb{P}^t$ generated by $p^t$. It is evident that both $\tau$ and $\tau^{-1}$ preserve measure, meaning that for any measurable set $A \subseteq \Sigma_m^+$,
$B \subseteq \Sigma_{m^t}^+$, $\mathbb{P}^t (\tau^{-1}A)=\mathbb{P} (A)$ and $\mathbb{P}^t (B)=\mathbb{P} (\tau B)$ hold. For any $\varpi=(j_1,j_2,\cdots)\in\Sigma_{m^t}^+$,
$j_s\in\{1,2,\cdots,m^t\}$, $s=1,2,\cdots$, we define
$$
  g_{\varpi,n}(x):= \begin{cases}
                    x, & n=0 \\
                    g_{j_n}\circ g_{j_{n-1}}\circ \cdots \circ g_{j_1}(x), & n\geq 1.
                   \end{cases}
$$
$(X, G^t)$ is referred to as the $t$-power system of $(X,G)$. Notably, for any $\varpi=(j_1,j_2,\cdots)\in\Sigma_{m^t}^+$,
$j_s\in\{1,2,\cdots,m^t\}$, $1\leq s$, there exists a unique $\omega=\tau(\varpi) \in \Sigma_m^+$, $\omega=(i_1,i_2,\cdots)$ such
that $f_{\omega, nt}(x)=g_{\varpi, n}(x)$ for any $n\geq 0$, any $x\in X$.

\subsection{Measure-theoretic entropy and topological entropy of free semigroup actions}
The measure-theoretic entropy and topological entropy of free semigroup actions have been extensively studied in the literature \cite{AR,CRV1,CRV2,K,KL,LMW,LW}. Here, we adopt the following definitions for measure-theoretic entropy and topological entropy.

\begin{definition}\cite{CRV1}
Let $G$ be a free semigroup acting on a compact metric space $X$, $\mu$ be G-invariant probability measure, $\xi$ be finite Borel measurable
	partition of X. Then the measure-theoretic entropy $h_{\mu}(G, \xi)$ of (X,G) with respect to $\xi$ is defined as
	$$
	h_{\mu}(G, \xi):=\lim_{k \rightarrow +\infty} \frac{1}{k} \int_{\Sigma_m^+} H_{\mu}(\bigvee_{i=0}^{k-1} f_{\omega, i}^{-1
	} \xi) d\mathbb{P}(\omega),
	$$
	and the measure-theoretic entropy of (X,G) is defined as
	$$
	h_{\mu}(G):=\sup_{\xi} h_{\mu}(G, \xi),
	$$
	where
	$$
	\begin{aligned}
		H_{\mu}(\bigvee_{i=0}^{k-1} f_{\omega, i}^{-1} \xi) & :=-\sum_{A \in \bigvee_{i=0}^{k-1} f_{\omega, i}^{-1} \xi} \mu(A) \log \mu
		(A) , \\
		\bigvee_{i=0}^{k-1} f_{\omega, i}^{-1} \xi & :=\xi \bigvee f_{\omega, 1}^{-1} \xi \bigvee \cdots \bigvee f_{\omega, k-1}^{-1} \xi.
	\end{aligned}
	$$
\end{definition}

We recall the concepts of separated sets and spanning sets in the context of free semigroup actions. Let $G$ be a free semigroup with $m$ generators $\{f_1,f_2,\cdots,f_m\}$ acting on compact metric space $X$, where $f_i$ is continuous self-map on $X$, $i \in \{1,2,\cdots,m\}$. A subset $E(\omega,k,\varepsilon) \subseteq X$ is defined as the largest cardinality $(\omega, k, \varepsilon)$ separated set if, for any distinct points $x,y \in E(\omega,k,\varepsilon)$, the distance $d^G_{\omega, k}(x,y) > \varepsilon$, and the cardinality $\sharp E(\omega,k,\varepsilon)$ of $E(\omega,k,\varepsilon)$ is maximized. Similarly, A subset $F(\omega,k,\varepsilon) \subseteq X$ is termed the smallest cardinality $(\omega, k, \varepsilon)$ spanning set if, for any $x \in X$, there exists $y \in F$ such that $d^G_{\omega, k}(x,y) \leq \varepsilon$, and the cardinality $\sharp F(\omega,k,\varepsilon)$ of $F(\omega,k,\varepsilon)$ is minimized.

\begin{definition}\cite{LD}
Let $G$ be a free semigroup acting on a compact metric space $X$. Then the topological entropy $h_{top}(G)$ of $(X,G)$ is
	$$
	h_{top}(G)=\lim_{\varepsilon \rightarrow 0} \liminf_{k \rightarrow +\infty} \frac{1}{k} \int_{\Sigma_m^+} \log \sharp E(\omega,k,\varepsilon) d\mathbb{P}(\omega).
	$$
\end{definition}

\begin{remark}\label{re2.1}
	Kifer \cite{K} introduced the topological entropy of random transformations as
	$$
	H_{top}(G)
	=\lim_{\varepsilon \rightarrow 0} \liminf_{k \rightarrow +\infty} \frac{1}{k}\log \sharp E(\omega,k,\varepsilon)
	$$
	for almost everywhere $\omega \in \Sigma_m^+$,
	and Li et al\cite{LD} demonstrated that $h_{top}(G)=H_{top}(G)$.
\end{remark}

Below, we present the ergodic theorem and the Brin-Katok local entropy formula for random dynamical systems, established by Kifer \cite{K} and Zhu \cite{Z1,Z2}, respectively
\begin{theorem}\label{ergodic}
\cite{K} Let $G$ be a free semigroup acting on a compact metric space $X$, $(\Sigma_m^+, \mathbb{P})$ be a probability space, $F: \Sigma_m^+ \times X \longrightarrow \Sigma_m^+ \times X $ be the skew
product transformation, $\mu$ be a Borel probability measure on $X$. If $\mathbb{P} \times \mu$ is ergodic with respect to $F$ and $\phi \in L^1(\mu)$, then there exists an $\Omega \subseteq \Sigma_m^+$ with $\mathbb{P}(\Omega)=1$, such that for any $\omega \in \Omega$, there exists a $X_{\omega} \subseteq X$ with $\mu(X_{\omega})=1$ where
$$
\lim_{k \rightarrow +\infty} \frac{1}{k} \sum_{i=0}^{k-1} \phi(f_{\omega,i}(x))=\int_X \phi(y) d \mu(y)
$$
holds for any $x \in X_{\omega}$.
\end{theorem}

\begin{remark}\label{re2.2}
It is also true that there exists a $W \subseteq X$ with $\mu(W)=1$, such that for any $ x \in W$, there exists an $\Omega_x \subseteq \Sigma_m^+$ with $\mathbb{P}(\Omega_x)=1$ where
$$
\lim_{k \rightarrow +\infty} \frac{1}{k} \sum_{i=0}^{k-1} \phi(f_{\omega,i}(x))=\int_X \phi(y) d \mu(y)
$$
holds for any $\omega \in \Omega_x$.
\end{remark}

\begin{theorem}\label{th2.3}
\cite{Z1,Z2} Let $G$ be a free semigroup acting on a compact metric space $X$, $(\Sigma_m^+, \mathbb{P})$ be a probability space, $F: \Sigma_m^+ \times X \longrightarrow \Sigma_m^+ \times X $ be the skew
product transformation, $\mu$ be a Borel probability measure on $X$. If $\mathbb{P} \times \mu$ is ergodic with respect to $F$ and $h_{\mu}(G)
< +\infty$, then there exists an $\Omega \subseteq \Sigma_m^+$ with $\mathbb{P}(\Omega)=1$ such that for any $\omega \in \Omega$, there
exists a $X_{\omega} \subseteq X$ with $\mu(X_{\omega})=1$ where
$$
\begin{aligned}
h_{\mu}(G)&=\lim_{\varepsilon \rightarrow 0} \liminf_{k \rightarrow +\infty} -\frac{1}{k} \log \mu (B_{\omega,k}^G(x,\varepsilon)) \\
          &=\lim_{\varepsilon \rightarrow 0} \limsup_{k \rightarrow +\infty} -\frac{1}{k} \log \mu (B_{\omega,k}^G(x,\varepsilon))
\end{aligned}
$$
holds for any $x \in X_{\omega}$.
\end{theorem}
\begin{remark}
The formulations of the two theorems differ in references \cite{K, Z1, Z2} due to modifications introduced for convenience. Readers are encouraged to independently verify the validity of these modifications.
\end{remark}

\section{Notions and properties}
In this Section, we introduce the concepts of correlation sum, upper(lower) local correlation entropy as well as $q-$order upper(lower) correlation entropy of free semigroup actions and explore their fundamental properties.\par
To begin, we introduce the concept of correlation sum for free semigroup actions.
\begin{definition}\label{def3.1}
	Let $G$ be a free semigroup acting on a compact metric space $X$. For any $x \in X$, $\varepsilon > 0$, $\omega \in \Sigma_m^+$, $k \geq
	1$, $n \geq 1$, the correlation sum of free semigroup actions is defined as follows
	$$
	\begin{aligned}
		C&(G,x,\varepsilon,\omega,k,n) \\
		&:=\frac{1}{n^2} \int_{\Sigma_m^+} \sharp \left\{(i,j): 0 \leq i, j \leq n-1, d_{\omega, k}^G(f_{\upsilon, i}(x), f_{\upsilon, j}(x))
		\leq \varepsilon \right\} d \mathbb{P} (\upsilon)
	\end{aligned}
	$$
	where $\sharp A$ is the cardinality of set $A$.
\end{definition}

\begin{remark}\label{re3.1}
\textcolor{red}{Initially, we propose an alternative definition to generalize the correlation sum, which is defined as}
$$
	\begin{aligned}
		C^{\prime}&(G,x,\varepsilon,\omega,k,n) \\
		&:=\frac{1}{n^2}  \sharp \left\{(i,j): 0 \leq i, j \leq n-1, d_{\omega, k}^G(f_{\omega, i}(x), f_{\omega, j}(x))
		\leq \varepsilon \right\}.
	\end{aligned}
$$
However, we prefer Definition \ref{def3.1} for the following reasons. Firstly, in classical dynamical systems, local correlation entropy is defined using a fixed Bowen metric $d_k$ to observe the first $n$ elements of the orbit $\{f^i(x)\}_{i=0}^{+\infty}$ and compute $C(x,\varepsilon, d_k, n)$. In the context of free semigroup actions, it is essential to choose a fixed Bowen metric, that is, $d_{\omega, k}^G$. Subsequently, this fixed Bowen metric is applied to observe the first n elements of the orbit $Orb\{x, G\}$, denoted as
$$
Orb(x,G,n):=\{f_{\omega,k}(x): \forall \omega \in \Sigma_m^+, 0 \leq k \leq n-1\}.
$$
We contend that $Orb\{x, G,n\}$ should encompass multiple trajectories induced by different $\upsilon$ rather than only unique trajectory induced by $\omega$. \textcolor{blue}{Secondly, if we adopt $C^{\prime}(G,x,\varepsilon,\omega,k,n) $ as the definition of generalized correlation sum,, we could not be able to obtain the analogue of Theorem (\ref{th1.1}) without additional conditions. The reason is that the ergodic theorem of random dynamical systems plays an important role in the proof of the Theorem (\ref{th1.1}).} Hence, the choice of the integral form is preferred.
\end{remark}

Building upon the concept of correlation sum outlined above, we introduce the local correlation entropy of free semigroup actions.
\begin{definition}
	Let $G$ be a free semigroup acting on a compact metric space $X$. The upper (lower) local correlation entropy of free semigroup actions is
	defined as follows
	$$
	\begin{aligned}
		&\overline{h}_{cor}(G,x):=\lim_{\varepsilon \rightarrow 0} \limsup_{k \rightarrow +\infty} - \frac{1}{k} \int_{\Sigma_m^+} \log
		\underline{C}(G,x,\varepsilon,\omega,k) d\mathbb{P}(\omega), \\
		&\underline{h}_{cor}(G,x):=\lim_{\varepsilon \rightarrow 0} \liminf_{k \rightarrow +\infty} - \frac{1}{k} \int_{\Sigma_m^+} \log
		\overline{C}(G,x,\varepsilon,\omega,k) d\mathbb{P}(\omega),
	\end{aligned}
	$$
	where
	$$
	\begin{aligned}
		&\underline{C}(G,x,\varepsilon,\omega,k):=\liminf_{n \rightarrow +\infty} C(G,x,\varepsilon,\omega,k,n),\\
		&\overline{C}(G,x,\varepsilon,\omega,k):= \limsup_{n \rightarrow +\infty} C(G,x,\varepsilon,\omega,k,n).
	\end{aligned}
	$$
	If $\overline{h}_{cor}(G,x) =\underline{h}_{cor}(G,x)$, then we denote $h_{cor}(G,x):=\overline{h}_{cor}(G,x) =\underline{h}_{cor}(G,x)$.
\end{definition}
Similarly, we introduce the concepts of correlation integral and correlation entropy of free semigroup actions as follows.
\begin{definition}
	Let $G$ be a free semigroup acting on a compact metric space $X$, $\mu$ be a Borel probability measure with full support on $X$. (this assumption holds when discussing correlation entropy). For $\varepsilon
	>0$, $k \geq 1$, and $q \in \mathbb{R}$, the correlation integral of $q$-order of free semigroup actions is defined as follows
	$$
	\begin{aligned}
		c(G,\mu,\varepsilon,k,q) &:= \frac{1}{q-1} \int_{\Sigma_m^+} \log \left(\int_X \mu(B_{\omega,k}^G(x,\varepsilon))^{q-1}
		d\mu(x)\right) d\mathbb{P}(\omega) \qquad & q \neq 1 ,\\
		c(G,\mu,\varepsilon,k,1) &:= \int_{\Sigma_m^+} \int_X \log \mu (B_{\omega, k}^G (x,\varepsilon)) d\mu(x) d\mathbb{P}(\omega) \qquad
		& q=1,
	\end{aligned}
	$$
	and the upper (lower) correlation entropy of $q$-order of free semigroup actions is defined as follows
	$$
	\begin{aligned}
		\overline{h}_{cor}(G,\mu,q) &: =\lim_{\varepsilon \rightarrow 0} \limsup_{k \rightarrow +\infty} -\frac{1}{k} c(G,\mu,\varepsilon,k,q),
		\\
		\underline{h}_{cor}(G,\mu,q) &: =\lim_{\varepsilon \rightarrow 0} \liminf_{k \rightarrow +\infty} -\frac{1}{k}
		c(G,\mu,\varepsilon,k,q).
	\end{aligned}
	$$
	If $\overline{h}_{cor}(G,\mu,q)=\underline{h}_{cor}(G,\mu,q)$, we denote
	$$
	h_{cor}(G,\mu,q):=\overline{h}_{cor}(G,\mu,q)=\underline{h}_{cor}(G,\mu,q).
	$$
\end{definition}

\begin{remark}
	When $G_{\ast}=\{f_1\}$, definition 2.1-2.3 degenerate into classical cases, as discussed in \cite{BTH,S2018,T}.
\end{remark}

\begin{remark}
Similar to the approach described in \cite{TV1998}, the definition of $q=1$ is imposed by continuity. For clarity, we provide an explanation here.
Let $\varepsilon >0$ satisfy the following condition, for any $\omega \in \Sigma_m^+$, $k \geq 1$ and $x \in X$,
$$
\mu\left( \{ y: d_{\omega,k}^G(x,y)=\varepsilon \} \right)=0.
$$
Assume that there exists $\omega_0 \in \Sigma_m^+$ such that $\inf_{x\in X}\mu(B_{\omega_0,k}^G(x,\varepsilon))=0$. Denote
$$
B_{\omega_0,k}^G(x,r_1,r_2):=\{y: r_1 < d_{\omega_0,k}^G(x,y) \leq r_2\}.
$$
\textcolor{blue}{It could be observed that
\begin{equation}\label{eq3.1}
\begin{aligned}
&\left|\mu(B_{\omega_0,k}^G(x,\varepsilon)) - \mu(B_{\omega_0,k}^G(y,\varepsilon))\right| \\
= & \left|\mu(B_{\omega_0,k}^G(x,\varepsilon) \setminus B_{\omega_0,k}^G(y,\varepsilon)) - \mu(B_{\omega_0,k}^G(y,\varepsilon) \setminus B_{\omega_0,k}^G(x,\varepsilon))\right| \\
\leq & \max \left\{ \mu(B_{\omega_0,k}^G(x,\varepsilon) \setminus B_{\omega_0,k}^G(y,\varepsilon)),  \mu(B_{\omega_0,k}^G(y,\varepsilon) \setminus B_{\omega_0,k}^G(x,\varepsilon)) \right\}.
\end{aligned}
\end{equation}
Note that if $d_{\omega_0,k}^G(x,y) < \delta < \varepsilon$, then we have
\begin{equation}\label{eq3.2}
\begin{aligned}
B_{\omega_0,k}^G(x,\varepsilon) \setminus B_{\omega_0,k}^G(y,\varepsilon) \subset B_{\omega_0,k}^G(x,\varepsilon - \delta, \varepsilon),\\
B_{\omega_0,k}^G(y,\varepsilon) \setminus B_{\omega_0,k}^G(x,\varepsilon) \subset B_{\omega_0,k}^G(x,\varepsilon, \varepsilon + \delta).
\end{aligned}
\end{equation}
Combined with the equations (\ref{eq3.1}) and (\ref{eq3.2}), we have
$$
\begin{aligned}
&\left|\mu(B_{\omega_0,k}^G(x,\varepsilon)) - \mu(B_{\omega_0,k}^G(y,\varepsilon))\right| \\
\leq & \max \left\{ \mu(B_{\omega_0,k}^G(x,\varepsilon) \setminus B_{\omega_0,k}^G(y,\varepsilon)),  \mu(B_{\omega_0,k}^G(y,\varepsilon) \setminus B_{\omega_0,k}^G(x,\varepsilon)) \right\} \\
\leq & \max \left\{ \mu(B_{\omega_0,k}^G(x,\varepsilon - \delta, \varepsilon)), \mu(B_{\omega_0,k}^G(x,\varepsilon, \varepsilon+ \delta)) \right\}.
\end{aligned}
$$
}
As $\delta$ approaches $0$, $B_{\omega_0,k}^G(x,\varepsilon , \varepsilon+ \delta)$ tends to $\emptyset$ and $B_{\omega_0,k}^G(x,\varepsilon - \delta, \varepsilon)$ tends to $\{ y: d_{\omega,k}^G(x,y)=\varepsilon \}$. Thus, for any $\eta >0$, there exists $\delta < \varepsilon$ such that if $d_{\omega_0,k}^G(x,y) < \delta$, then
$$
\max \left\{ \mu(B_{\omega_0,k}^G(x,\varepsilon - \delta, \varepsilon)) \mu(B_{\omega_0,k}^G(x,\varepsilon, \varepsilon+ \delta)) \right\}\leq \eta.
$$
Hence, $\mu(B_{\omega_0,k}^G(x,\varepsilon))$ is continuous with respect to $x$. Moreover, $X$ is compact, implying the existence of $x_0 \in X$ such that $\mu(B_{\omega_0,k}^G(x_0,\varepsilon))=0$. However, this contradicts the condition that the support of $\mu$ is $X$. Therefore, $0<\inf_{x \in X}\mu(B_{\omega,k}^G(x,\varepsilon)) \leq 1$ for any $\omega \in \Sigma_m^+$, ensuring the interchangeability of the limit and the integral in the following process
$$
\begin{aligned}
&\lim_{q \rightarrow 1}\frac{1}{q-1} \int_{\Sigma_m^+} \log \left(\int_X \mu(B_{\omega,k}^G(x,\varepsilon))^{q-1} d\mu(x)\right)
d\mathbb{P}(\omega)
\\ =&\int_{\Sigma_m^+}\lim_{q \rightarrow 1}\frac{1}{q-1} \log \left(\int_X \mu(B_{\omega,k}^G(x,\varepsilon))^{q-1} d\mu(x)\right)
d\mathbb{P}(\omega) \\
=&\int_{\Sigma_m^+}\lim_{q \rightarrow 1}\frac{\int_X \mu(B_{\omega,k}^G(x,\varepsilon))^{q-1}\log \mu (B_{\omega, k}^G (x,\varepsilon)) d\mu(x)}{\int_X \mu(B_{\omega,k}^G(x,\varepsilon))^{q-1} d\mu(x)}d\mathbb{P}(\omega)\\
=&\int_{\Sigma_m^+} \int_X \log \mu (B_{\omega, k}^G (x,\varepsilon)) d\mu(x) d\mathbb{P}(\omega).
\end{aligned}
$$
Furthermore, we will demonstrate in the proof of Theorem \ref{th1.1} that the set of such $\varepsilon$ forms an uncountable dense subset of $\mathbb{R}$.
\end{remark}
\begin{proposition}\label{prop3.1}
	(1) If $ 0 < \varepsilon_1 \leq \varepsilon_2$, then
	$C(G,x,\varepsilon_1,\omega,k,n) \leq C(G,x,\varepsilon_2,\omega,k,n).
	$
	So $\overline{h}_{cor}(G,x)$ is well defined. \\
	(2) If $k_1 \leq k_2$, then
	$$
	\int_{\Sigma_m^+} \log \liminf_{n \rightarrow +\infty} C(G,x,\varepsilon,\omega,k_1,n) d\mathbb{P}(\omega) \geq \int_{\Sigma_m^+}
	\log \liminf_{n \rightarrow +\infty} C(G,x,\varepsilon,\omega,k_2,n) d\mathbb{P}(\omega).
	$$
	(3) If $\{n_i\}_i$ increases strictly and $\frac{n_{i+1}}{n_i} \rightarrow 1$, then
	$$
	\liminf_{n \rightarrow +\infty} C(G,x,\varepsilon,\omega,k,n)=\liminf_{i \rightarrow +\infty} C(G,x,\varepsilon,\omega,k,n_i).
	$$
	(4) If $\{k_i\}_i$ increases strictly and $\frac{k_{i+1}}{k_i} \rightarrow 1$, then
	$$
	\overline{h}_{cor}(G,x)=\lim_{\varepsilon \rightarrow 0} \limsup_{i \rightarrow +\infty} - \frac{1}{k_i} \int_{\Sigma_m^+} \log \liminf_{n
		\rightarrow +\infty} C(G,x,\varepsilon,\omega,k_i,n) d\mathbb{P}(\omega).
	$$
	\begin{proof}
		(1) and (2) can be proved by definition. \\
		(3) For any given $n$, there exists an index $s$ such that $n_s \leq n \leq n_{s+1}$. Thus, we have
		$$
		\begin{aligned}
			\frac{n_s^2}{n_{s+1}^2} & C(G,x,\varepsilon,\omega,k,n_s) \\
			& \leq \frac{1}{n^2} \int_{\Sigma_m^+} \sharp \{(i,j): 0 \leq i, j \leq n_s-1, d_{\omega, k}^G(f_{\upsilon, i}(x), f_{\upsilon,
				j}(x)) \leq \varepsilon \} d \mathbb{P} (\upsilon) \\
			& \leq C(G,x,\varepsilon,\omega,k,n) \\
			& \leq \frac{1}{n^2} \int_{\Sigma_m^+} \sharp \{(i,j): 0 \leq i, j \leq n_{s+1}-1, d_{\omega, k}^G(f_{\upsilon, i}(x), f_{\upsilon,
				j}(x)) \leq \varepsilon \} d \mathbb{P} (\upsilon) \\
			& \leq \frac{n_{s+1}^2}{n_s^2} C(G,x,\varepsilon,\omega,k,n_{s+1}).
		\end{aligned}
		$$
		Taking the $\liminf$ on both sides yields the desired result. \\
		(4) By employing (2) and following the methodology outlined in (3), the proof is established.
	\end{proof}
	\begin{remark}
It extends the results of V.$\check{S}$pitalsk$\check{y}$ \cite{S2018}. All the properties remain valid when considering $\overline{C}(G,x,\varepsilon,\omega,k)$ and $\underline{h}_{cor}(G,x)$.
	\end{remark}
\end{proposition}
\begin{proposition}\label{prop3.2}
	(1) If $0 < \varepsilon_1 \leq \varepsilon_2$, then $c(G,\mu,\varepsilon_1,k,q) \leq c(G,\mu,\varepsilon_2,k,q)$. So
	$\overline{h}_{cor}(G,\mu,q)$ is well-defined. \\
	(2) If $k_1 \leq k_2$, then $c(G,\mu,\varepsilon,k_1,q) \geq c(G,\mu,\varepsilon,k_2,q)$. \\
	(3) If $\{k_i\}_i$ increases strictly and satisfies $\frac{k_{i+1}}{k_i} \rightarrow 1$, then \\
	$$
	\overline{h}_{cor}(G,\mu,q)=\lim_{\varepsilon \rightarrow 0} \limsup_{i \rightarrow +\infty} -\frac{1}{k_i}c(G,\mu,\varepsilon,k_i,q).
	$$
	(4) If $q_1 < q_2$, then $c(G,\mu,\varepsilon,k,q_1) \leq c(G,\mu,\varepsilon,k,q_2)$. In particularly,
	$$
		\overline{h}_{cor}(G,\mu,q_1) \geq \overline{h}_{cor}(G,\mu,q_2) ,\quad
	    \underline{h}_{cor}(G,\mu,q_1) \geq \underline{h}_{cor}(G,\mu,q_2).
	$$
	\begin{proof}
		(1) If $q > 1$, then
		$$
		\begin{aligned}
			&\int_{\Sigma_m^+} \log \left(\int_X \mu(B_{\omega,k}^G(x,\varepsilon_1))^{q-1} d\mu(x)\right) d\mathbb{P}(\omega) \\
			\leq &\int_{\Sigma_m^+} \log \left(\int_X \mu(B_{\omega,k}^G(x,\varepsilon_2))^{q-1} d\mu(x)\right) d\mathbb{P}(\omega).
		\end{aligned}
		$$
		So $c(G,\mu,\varepsilon_1,k,q) \leq c(G,\mu,\varepsilon_2,k,q)$. If $q=1$, then
		$$
		\int_{\Sigma_m^+} \int_X \log \mu (B_{\omega, k}^G (x,\varepsilon_1)) d\mu(x) d\mathbb{P}(\omega) \leq \int_{\Sigma_m^+} \int_X \log
		\mu (B_{\omega, k}^G (x,\varepsilon_2)) d\mu(x) d\mathbb{P}(\omega).
		$$
		So $c(G,\mu,\varepsilon_1,k,q) \leq c(G,\mu,\varepsilon_2,k,q)$. If $q < 1$, then
		$$
		\begin{aligned}
			&\int_{\Sigma_m^+} \log \left(\int_X \mu(B_{\omega,k}^G(x,\varepsilon_1))^{q-1} d\mu(x)\right) d\mathbb{P}(\omega) \\
			\geq &\int_{\Sigma_m^+} \log \left(\int_X \mu(B_{\omega,k}^G(x,\varepsilon_2))^{q-1} d\mu(x)\right) d\mathbb{P}(\omega).
		\end{aligned}
		$$
		So $c(G,\mu,\varepsilon_1,k,q) \leq c(G,\mu,\varepsilon_2,k,q)$. \\
		(2) It can be proved by definition.\\
		(3) It can be proved in a similar manner as Proposition \ref{prop3.1}. \\
		(4) Given that $\mu$ possesses full support, for $q <1$, the inequality
		$$
		(q-1) \int_X \log \mu (B_{\omega, k}^G(x,\varepsilon))d\mu(x) \leq \log \int_X \mu (B_{\omega,k}^G(x,\varepsilon))^{q-1}d\mu(x)
		$$
		follows from Jensen's inequality \cite{S1984}. Consequently,
$$
c(G,\mu,\varepsilon,k,1) > c(G,\mu,\varepsilon,k,q).
$$
Similarly, for $q>1$, we obtain $c(G,\mu,\varepsilon,k,1) < c(G,\mu,\varepsilon,k,q)$. Furthermore, for $ 1 < q_1< q_2 $, the inequality
		$$
		\begin{aligned}
			\frac{1}{q_1-1} \int_{\Sigma_m^+} \log \left(\int_X \mu (B_{\omega,k}^G(x,\varepsilon))^{q_1-1}d\mu(x)\right) d\mathbb{P}(\omega) \\
			< \frac{1}{q_2-1} \int_{\Sigma_m^+} \log \left(\int_X \mu (B_{\omega,k}^G(x,\varepsilon))^{q_2-1}d\mu(x)\right) d\mathbb{P}(\omega)
		\end{aligned}
		$$
		holds by Lyapunov's inequality \cite{S1984}. Since $\mu$ has full support, Lyapunov's inequality applies to $q_1<q_2<1$, thus concluding the proof.
	\end{proof}
	\begin{remark}
		It generalizes the results of E.Verbitskiy \cite{V} and all the properties remain valid when considering $\underline{h}_{cor}(G,\mu,q)$.
	\end{remark}
\end{proposition}
\begin{proposition}
	For any integer $t\geq 1$, let $G$ be a free semigroup acting on a compact metric space $X$, $(X,G^t)$ $t$-power system of $(X,G)$. Then for any $q \in \mathbb{R}$
	$$
		\overline{h}_{cor}(G^t,\mu,q)  = t \cdot \overline{h}_{cor}(G,\mu,q),\quad
		\underline{h}_{cor}(G^t,\mu,q)  =t \cdot \underline{h}_{cor}(G,\mu,q).
	$$
	\begin{proof}
		Since $X$ is a compact metric space and each $f_i \in G_{\ast}$ is continuous, for any $\varepsilon > 0$, there exists $\delta \leq
		\varepsilon$ such that $d(x,y) \leq \delta$ implies $d_{\omega ,t}^G(x,y) \leq \varepsilon$ for any $\omega \in \Sigma_m^+$.
		Therefore, for any $\varpi \in \Sigma_{m^t}^+$, there exists a unique $\omega \in \Sigma_m^+$ such that
		$$
		\mu (B_{\omega,tk}^G(x,\delta)) \leq \mu (B_{\varpi,k}^{G^t}(x,\delta)) \leq \mu (B_{\omega,tk}^G(x,\varepsilon)).
		$$
		By computation, we establish for any $q \in \mathbb{R}$,
$$
c(G,\mu,\delta,tk,q) \leq c(G^t,\mu,\delta,k,q)\leq c(G,\mu,\varepsilon,tk,q).
$$
Combined
		with Proposition \ref{prop3.2}, we derive
		$$
			\overline{h}_{cor}(G^t,\mu,q)  = t \cdot \overline{h}_{cor}(G,\mu,q),\quad
			\underline{h}_{cor}(G^t,\mu,q)  =t \cdot \underline{h}_{cor}(G,\mu,q).
		$$
	\end{proof}
\end{proposition}

\section{Proof of Theorem \ref{th1.1}}
In this Section, we present the proof of Theorem \ref{th1.1} ($q=2$) following the methodology from \cite{MS}. The proof of Theorem \ref{th1.1} proceeds in three main steps. \textcolor{red}{Firstly, we establish the existence of a countable subset $Q \subset \mathbb{R}$ such that for any $\varepsilon \notin Q$, $\omega \in \Sigma_m^+$ and $k \in \mathbb{N}$,}
$$
\mu \times \mu \left(  \left\{ (x,y)  \in X \times X: d_{\omega, k}^G(x,y)= \varepsilon \right\} \right)=0,
$$
meaning that the mapping $\varepsilon \mapsto \int_X \mu ( B_{\omega, k}^G(x, \varepsilon) ) d\mu(x)$ is continuous at $\varepsilon$. Secondly, for a fixed $\omega \in \Sigma_m^+$, we identify a full measure subset $W(\omega, k, \varepsilon) \subseteq X$ where the equality (\ref{eq1.2}) is established. Finally, we determine a common full measure subset applicable for any $\omega \in \Sigma_m^+$. We commence by proving the first step.\par
\begin{proof}[Proof of Theorem \ref{th1.1}]
\textbf{Step 1}. \textcolor{red}{For any $\omega \in \Sigma_m^+$ and $k \in \mathbb{N}$, the set of real number $\varepsilon \in \mathbb{R}$ satisfying the inequality}
$$
\mu \times \mu \left(  \left\{ (x,y)  \in X \times X: d_{\omega, k}^G(x,y)= \varepsilon \right\} \right)> \frac{1}{n}
$$
has a cardinality less than $n$, owing to $\mu \times \mu$ being a probability measure on the compact metric space $X \times X$. Let $Q_{\omega, k, n}$ denote the collection of such $\varepsilon$. Consequently, there exists a countable subset $\bigcup_{n=1}^{+\infty} Q_{\omega, k, n}$ such that for any $\varepsilon \notin \bigcup_{n=1}^{+\infty} Q_{\omega, k, n}$,
$$
\mu \times \mu \left(  \left\{ (x,y)  \in X \times X: d_{\omega, k}^G(x,y)= \varepsilon \right\} \right)=0.
$$
Given that $\Sigma_m^+$ is a compact metric space, there exists a countable dense subset $\{ \omega_r \}_{r=1}^{+\infty} \subset \Sigma_m^+$. For any $\omega \in \Sigma_m^+$ and $k \in \mathbb{N}$, there exists an $\omega_{r_0}$ such that the Bowen metric $d_{\omega, k}^G$ equals the Bowen metric $d_{\omega_{r_0}, k}^G$, implying that
$$
\bigcup_{n=1}^{+\infty} Q_{\omega, k, n} = \bigcup_{n=1}^{+\infty} Q_{\omega_{r_0}, k, n}.
$$
Thus, there exists a countable subset $Q:=\bigcup_{r,k,n=1}Q_{\omega_r, k, n}$ such that for any $\varepsilon \notin Q$, $\omega \in \Sigma_m^+$ and $k \in \mathbb{N}$,
$$
\mu \times \mu \left(  \left\{ (x,y)  \in X \times X: d_{\omega, k}^G(x,y)= \varepsilon \right\} \right)=0.
$$
It is noted that
$$
\int_X \mu(B_{\omega, k}^G(x,\varepsilon))d \mu(x) = \mu \times \mu \left(  \left\{ (x,y)  \in X \times X: d_{\omega, k}^G(x,y)\leq \varepsilon \right\} \right).
$$
For any sequence $\{\varepsilon_n\}_{n=1}^{+\infty}$ with $\varepsilon_n < \varepsilon$ and $\lim_{n \rightarrow +\infty}\varepsilon_n =\varepsilon$, we have
$$
\begin{aligned}
&\lim_{n \rightarrow +\infty}\left(\int_X \mu(B_{\omega, k}^G(x,\varepsilon))d \mu(x)-\int_X \mu(B_{\omega, k}^G(x,\varepsilon_n))d \mu(x)\right)\\
=&\lim_{n \rightarrow +\infty}\mu \times \mu \left(  \left\{ (x,y)  \in X \times X: \varepsilon_n < d_{\omega, k}^G(x,y)\leq \varepsilon \right\} \right)\\
=&\mu \times \mu \left(  \left\{ (x,y)  \in X \times X: d_{\omega, k}^G(x,y)= \varepsilon \right\} \right)\\
=&0.
\end{aligned}
$$
Similarly, for any sequence $\{\varepsilon_n\}_{n=1}^{+\infty}$ with $\varepsilon_n > \varepsilon$ and $\lim_{n \rightarrow +\infty}\varepsilon_n =\varepsilon$, we obtain
$$
\begin{aligned}
&\lim_{n \rightarrow +\infty}\left(\int_X \mu(B_{\omega, k}^G(x,\varepsilon_n))d \mu(x)-\int_X \mu(B_{\omega, k}^G(x,\varepsilon))d \mu(x)\right)\\
=&\lim_{n \rightarrow +\infty}\mu \times \mu \left(  \left\{ (x,y)  \in X \times X: \varepsilon < d_{\omega, k}^G(x,y)\leq \varepsilon_n \right\} \right)\\
=&0.
\end{aligned}
$$
Hence, for any $\omega \in \Sigma_m^+$ and $k \in \mathbb{N}$, the mapping $\varepsilon \mapsto \int_X \mu ( B_{\omega, k}^G(x, \varepsilon) ) d\mu(x)$ is continuous at $\varepsilon \notin Q$.

\textbf{Step 2}. We claim that for any given $\varepsilon \notin Q$, $\omega \in \Sigma_m^+$ and $k\in \mathbb{N}$, there exists a subset $W(\omega, k, \varepsilon) \subseteq X$ of full measure such that for any $x \in W(\omega, k, \varepsilon)$, the following convergence holds
$$
\lim_{n \rightarrow +\infty} C(G,x,\varepsilon,\omega,k,n) = \int_X \mu (B_{\omega,k}^G(x,\varepsilon))d\mu(x).
$$
\par
\begin{proof}
Given $\omega \in \Sigma_m^+$, $k\in \mathbb{N}$ and $\varepsilon \notin Q$, for any $t \in \mathbb{N}$, there exists a finite measurable partition of $X$,
denoted by $\xi_t:=\{A_{t,1}, A_{t,2}, \cdots, A_{t,N(t)}\}$, satisfying
$$
\mu(A_{t,1}) \leq 2^{-t}, \quad {\rm diam}_{\omega,k}(A_{t,s}) \leq 2^{-t}, \quad 2 \leq s \leq N(t),
$$
where ${\rm diam}_{\omega,k}(A):=\sup_{x, y \in A}d_{\omega,k}^G(x,y)$. Furthermore, we could consider that the boundary of $A_{t,s}$, $1 \leq s \leq N(t)$, to have measure $0$, given that $\mu$ is Borel probability measure and $X$ is compact metric space.\par
Denote
$$
\begin{aligned}
S_{\varepsilon}&:=\{(x,y) \in X \times X | d_{\omega,k}^G(x,y) \leq \varepsilon\}, \\
\mathcal{C}_1&:=\{A_{t,s_1} \times A_{t,s_2} \in \xi_t \times \xi_t | A_{t,s_1} \times A_{t,s_2} \subseteq S_{\varepsilon} \}, \\
\mathcal{C}_2&:=\{A_{t,s_1} \times A_{t,s_2} \in \xi_t \times \xi_t | A_{t,s_1} \times A_{t,s_2} \cap S_{\varepsilon} \neq \emptyset
\}.
\end{aligned}
$$
It is obvious that
$$
\bigcup_{A_{t,s_1} \times A_{t,s_2} \in \mathcal{C}_1} A_{t,s_1} \times A_{t,s_2} \subseteq S_{\varepsilon} \subseteq
\bigcup_{A_{t,s_1} \times A_{t,s_2} \in \mathcal{C}_2} A_{t,s_1} \times A_{t,s_2}.
$$
In particularly, we claim
\begin{equation}\label{eq4.1}
\begin{aligned}
S_{\varepsilon-2^{-t+1}} & - ((A_{t,1} \times X) \cup (X \times A_{t,1})) \\
& \subseteq \bigcup_{A_{t,s_1} \times A_{t,s_2} \in \mathcal{C}_1} A_{t,s_1} \times A_{t,s_2} \\
& \subseteq \bigcup_{A_{t,s_1} \times A_{t,s_2} \in \mathcal{C}_2} A_{t,s_1} \times A_{t,s_2} \\
& \subseteq S_{\varepsilon+2^{-t+1}} \cup (A_{t,1} \times X) \cup (X \times A_{t,1}).
\end{aligned}
\end{equation}
Indeed, for any $(x,y) \in S_{\varepsilon-2^{-t+1}} - ((A_{t,1} \times X) \cup (X \times A_{t,1})),$ it holds that
$$
d_{\omega,k}^G(x,y) \leq \varepsilon-2^{-t+1}, \quad x \notin A_{t,1}, \quad y \notin A_{t,1}.
$$
Hence, there exists $s_1,s_2 \neq 1$ such that $(x,y) \in A_{t,s_1} \times A_{t,s_2}$. For any $ (x^{\prime},y^{\prime}) \in
A_{t,s_1} \times A_{t,s_2} $, we have
$$
d_{\omega,k}^G(x^{\prime},y^{\prime}) \leq d_{\omega,k}^G(x^{\prime},x) + d_{\omega,k}^G(x,y) +d_{\omega,k}^G(y,y^{\prime}) \leq
2^{-t} + \varepsilon-2^{-t+1} + 2^{-t} = \varepsilon ,
$$
meaning that $(x^{\prime},y^{\prime}) \in S_{\varepsilon}$. Therefore, $A_{t,s_1} \times A_{t,s_2} \in \mathcal{C}_1$, which yields
$$
S_{\varepsilon-2^{-t+1}} - ((A_{t,1} \times X) \cup (X \times A_{t,1})) \subseteq \bigcup_{A_{t,s_1} \times A_{t,s_2} \in
\mathcal{C}_1} A_{t,s_1} \times A_{t,s_2} .
$$
Similarly, it can be demonstrated that
$$
\bigcup_{A_{t,s_1} \times A_{t,s_2} \in \mathcal{C}_2} A_{t,s_1} \times A_{t,s_2}  \subseteq S_{\varepsilon+2^{-t+1}} \cup (A_{t,1} \times X) \cup (X \times A_{t,1}).
$$
According to Theorem \ref{ergodic}, for any characteristic function $\chi_{_{A_{t,s}}}$, where $1 \leq s \leq N(t)$, there exists a subset $W_{t,s} \subseteq X$ with $\mu(W_{t,s})=1$ such that for any $x \in W_{t,s}$, there exists a subset $\Omega_{x,t,s} \subseteq \Sigma_m^+$ with $\mathbb{P}(\Omega_{x,t,s})=1$ satisfying for any $\upsilon \in \Omega_{x,t,s}$,
$$
\lim_{n \rightarrow +\infty} \frac{1}{n} \sum_{i=0}^{n-1} \chi_{_{A_{t,s}}}(f_{\upsilon,i}(x))=\mu(A_{t,s}).
$$
Furthermore, utilizing Egoroff's theorem\cite{H}, for any $\delta > 0$, there exists a subset $\Omega_{\delta,x,t,s} \subseteq \Omega_{x,t,s}$ with $\mathbb{P}(\Omega_{\delta,x,t,s}) > 1 - \frac{\delta}{N(t)}$. This subset $\Omega_{\delta,x,t,s}$ satisfies the existence of $N(\Omega_{\delta,x,t,s}) \in \mathbb{N}$ such that if $n > N(\Omega_{\delta,x,t,s})$, then for any $\upsilon \in \Omega_{\delta,x,t,s}$, the following inequality holds
$$
\left|\frac{1}{n} \sum_{i=0}^{n-1} \chi_{_{A_{t,s}}}(f_{\upsilon,i}(x)) - \mu(A_{t,s})\right| \leq \frac{2^{-t-1}}{N^2(t)}.
$$
Let $W_t:=\bigcap_{s=1}^{N(t)} W_{t,s}$, where it is evident that $\mu(W_t)=1$. For any $x \in W_t$, define
$
\Omega_{\delta,x,t} := \bigcap_{s=1}^{N(t)} \Omega_{\delta,x,t,s}.
$
It is noteworthy that $\mathbb{P}(\Omega_{\delta,x,t}) \geq 1 - \delta$. Define $$N(\Omega_{\delta,x,t}):=\max_{1 \leq s \leq N(t)} \{ N(\Omega_{\delta,x,t,s})\}.$$ Hence, if $n > N(\Omega_{\delta,x,t})$, then for any $\upsilon \in \Omega_{\delta,x,t}$ and $1 \leq s \leq N(t)$,
$$
\left|\frac{1}{n} \sum_{i=0}^{n-1} \chi_{_{A_{t,s}}}(f_{\upsilon,i}(x)) - \mu(A_{t,s})\right| \leq \frac{2^{-t-1}}{N^2(t)}.
$$
Consequently, for any $x \in W_t$, $\delta > 0$, if $n > N(\Omega_{\delta,x,t})$, then for any $\upsilon \in \Omega_{\delta,x,t}$ and $1 \leq s_1, s_2 \leq N(t)$, we have
\begin{equation}\label{4.2}
\begin{aligned}
& \left|\frac{1}{n^2} \sharp \{(i,j) : 0 \leq i,j \leq n-1, (f_{\upsilon,i}(x),f_{\upsilon,j}(x)) \in A_{t,s_1} \times A_{t,s_2} \} \right. \\
& - \mu \times \mu (A_{t,s_1} \times A_{t,s_2}) \bigg| \\
= & \left| \frac{1}{n} \sharp \{i : 0 \leq i \leq n-1, f_{\upsilon,i}(x) \in A_{t,s_1} \} \cdot \frac{1}{n} \sharp \{j : 0 \leq j \leq n-1, f_{\upsilon,j}(x) \in A_{t,s_2} \} \right. \\
&- \mu(A_{t,s_1}) \cdot \mu(A_{t,s_2}) \bigg| \\
\leq & \frac{2^{-t}}{N^2(t)}.
\end{aligned}
\end{equation}
Note that
\begin{equation}\label{4.3}
\begin{aligned}
& \frac{1}{n^2} \sharp \left( \bigcup_{A_{t,s_1} \times A_{t,s_2} \in \mathcal{C}_1} \{(i,j) | 0 \leq i,j \leq n-1,
(f_{\upsilon,i}(x),f_{\upsilon,j}(x)) \in A_{t,s_1} \times A_{t,s_2} \} \right) \\
\leq & \frac{1}{n^2} \sharp \{(i,j) | 0 \leq i,j \leq n-1, d_{\omega,k}^G(f_{\upsilon,i}(x),f_{\upsilon,j}(x)) \leq \varepsilon \} \\
\leq & \frac{1}{n^2} \sharp \left( \bigcup_{A_{t,s_1} \times A_{t,s_2} \in \mathcal{C}_2} \{(i,j) | 0 \leq i,j \leq n-1,
(f_{\upsilon,i}(x),f_{\upsilon,j}(x)) \in A_{t,s_1} \times A_{t,s_2} \} \right).
\end{aligned}
\end{equation}
In light of (\ref{eq4.1}), we derive
\begin{equation}\label{4.4}
\begin{aligned}
(\mu\times\mu) \left( S_{\varepsilon-2^{-t+1}} \right) - 2 \mu (A_{t,1}) &\leq \sum_{A_{t,s_1} \times A_{t,s_2} \in \mathcal{C}_1} (\mu\times\mu) \left(A_{t,s_1} \times A_{t,s_2} \right) \\
&\leq \sum_{A_{t,s_1} \times A_{t,s_2} \in \mathcal{C}_2} (\mu\times\mu) \left(A_{t,s_1} \times A_{t,s_2} \right)\\
&\leq (\mu\times\mu) \left( S_{\varepsilon+2^{-t+1}} \right) +2\mu \left( A_{t,1} \right).
\end{aligned}
\end{equation}
By combining formulas (\ref{4.2}), (\ref{4.3}), and (\ref{4.4}), we derive the following inequality
$$
\begin{aligned}
& \frac{1}{n^2} \sharp \{(i,j) | 0 \leq i,j \leq n-1, d_{\omega,k}^G(f_{\upsilon,i}(x),f_{\upsilon,j}(x)) \leq \varepsilon \} \\
\geq & \sum_{A_{t,s_1} \times A_{t,s_2} \in \mathcal{C}_1} \frac{1}{n^2} \sharp \{(i,j) : 0 \leq i,j \leq n-1,
(f_{\upsilon,i}(x),f_{\upsilon,j}(x)) \in A_{t,s_1} \times A_{t,s_2} \} \\
\geq & \sum_{A_{t,s_1} \times A_{t,s_2} \in \mathcal{C}_1} \left( (\mu\times\mu) \left(A_{t,s_1} \times A_{t,s_2} \right) - \frac{2^{-t}}{N^2(t)} \right)\\
= & \sum_{A_{t,s_1} \times A_{t,s_2} \in \mathcal{C}_1} (\mu\times\mu) \left(A_{t,s_1} \times A_{t,s_2} \right) - \sharp \mathcal{C}_1 \frac{2^{-t}}{N^2(t)} \\
\geq & (\mu\times\mu)\left( S_{\varepsilon-2^{-t+1}} \right) -2 \mu (A_{t,1}) -2^{-t}\\
\geq & (\mu\times\mu)\left( S_{\varepsilon-2^{-t+1}} \right)-3\cdot2^{-t}.
\end{aligned}
$$
Similarly, we could get the another inequality as follows,
$$
\begin{aligned}
& (\mu\times\mu) (S_{\varepsilon-2^{-t+1}}) - 3 \cdot 2^{-t} \\
\leq & \frac{1}{n^2} \sharp \{(i,j) | 0 \leq i,j \leq n-1, d_{\omega,k}^G(f_{\upsilon,i}(x),f_{\upsilon,j}(x)) \leq \varepsilon \} \\
\leq & (\mu\times\mu) (S_{\varepsilon+2^{-t+1}}) + 3 \cdot 2^{-t}.
\end{aligned}
$$
Recalling the definition of correlation sum under free semigroup actions, we have
$$
\begin{aligned}
	C&(G,x,\varepsilon,\omega,k,n) \\
	=&\frac{1}{n^2} \int_{\Sigma_m^+} \sharp \{(i,j): 0 \leq i, j \leq n-1, d_{\omega, k}^G(f_{\upsilon, i}(x), f_{\upsilon, j}(x)) \leq \varepsilon \} d \mathbb{P} (\upsilon)\\
	=&\frac{1}{n^2} \int_{\Omega_{\delta,x,t}} \sharp \{(i,j): 0 \leq i, j \leq n-1, d_{\omega, k}^G(f_{\upsilon, i}(x), f_{\upsilon, j}(x)) \leq \varepsilon \} d \mathbb{P} (\upsilon) \\
	&+\frac{1}{n^2} \int_{\Sigma_m^+ - \Omega_{\delta,x,t}} \sharp \{(i,j): 0 \leq i, j \leq n-1, d_{\omega, k}^G(f_{\upsilon, i}(x), f_{\upsilon, j}(x)) \leq \varepsilon \} d \mathbb{P} (\upsilon).
\end{aligned}
$$
Now, we provide an estimation of the correlation sum. For any $x \in W_t$, any $\delta > 0$, if $n > N(\Omega_{\delta, x,t})$, then
$$
\begin{aligned}
& \int_{\Omega_{\delta,x,t}} (\mu\times\mu) (S_{\varepsilon-2^{-t+1}}) - 3 \cdot 2^{-t} d \mathbb{P} (\upsilon) \\
\leq & \frac{1}{n^2} \int_{\Omega_{\delta,x,t}} \sharp \{(i,j): 0 \leq i, j \leq n-1, d_{\omega, k}^G(f_{\upsilon, i}(x), f_{\upsilon, j}(x)) \leq \varepsilon \} d \mathbb{P} (\upsilon) \\
\leq & \int_{\Omega_{\delta,x,t}} (\mu\times\mu) (S_{\varepsilon+2^{-t+1}}) + 3 \cdot 2^{-t} d \mathbb{P} (\upsilon),
\end{aligned}
$$
which implies
$$
\begin{aligned}
& \left( (\mu\times\mu) (S_{\varepsilon-2^{-t+1}}) - 3 \cdot 2^{-t} \right) \mathbb{P} (\Omega_{\delta,x,t}) \\
\leq & \frac{1}{n^2} \int_{\Omega_{\delta,x,t}} \sharp \{(i,j): 0 \leq i, j \leq n-1, d_{\omega, k}^G(f_{\upsilon, i}(x), f_{\upsilon, j}(x)) \leq \varepsilon \} d \mathbb{P} (\upsilon) \\
\leq & \left( (\mu\times\mu) (S_{\varepsilon+2^{-t+1}}) + 3 \cdot 2^{-t} \right) \mathbb{P} (\Omega_{\delta,x,t}).
\end{aligned}
$$
Therefore,
$$
\begin{aligned}
	& \left( (\mu\times\mu) (S_{\varepsilon-2^{-t+1}}) - 3 \cdot 2^{-t} \right) \mathbb{P} (\Omega_{\delta,x,t}) \\
	\leq & C(G,x,\varepsilon,\omega,k,n) \\
	\leq & \left( (\mu\times\mu) (S_{\varepsilon+2^{-t+1}}) + 3 \cdot 2^{-t} \right) \mathbb{P} (\Omega_{\delta,x,t}) \\
	& + \frac{1}{n^2} \int_{\Sigma_m^+ - \Omega_{\delta,x,t}} \sharp \{(i,j): 0 \leq i, j \leq n-1, d_{\omega, k}^G(f_{\upsilon, i}(x), f_{\upsilon, j}(x)) \leq \varepsilon \} d \mathbb{P} (\upsilon) \\
	\leq & \left( (\mu\times\mu) (S_{\varepsilon+2^{-t+1}}) + 3 \cdot 2^{-t} \right) \mathbb{P} (\Omega_{\delta,x,t}) + \mathbb{P} (\Sigma_m^+ - \Omega_{\delta,x,t}).
\end{aligned}
$$
Given that $ \mathbb{P} ( \Omega_{\delta,x,t} ) \geq 1 - \delta$ and $ \mathbb{P} (\Sigma_m^+ - \Omega_{\delta,x,t}) \leq \delta$, we can establish the following inequlity
$$
\begin{aligned}
& \left( (\mu\times\mu) (S_{\varepsilon-2^{-t+1}}) - 3 \cdot 2^{-t} \right) \left( 1 - \delta \right) \\
\leq & C(G,x,\varepsilon,\omega,k,n) \\
\leq & \left( (\mu\times\mu) (S_{\varepsilon+2^{-t+1}}) + 3 \cdot 2^{-t} \right) + \delta.
\end{aligned}
$$
Taking $\liminf$ and $\limsup$ of $C(G,x,\varepsilon,\omega,k,n)$ as $n \rightarrow + \infty$, we obtain
$$
\begin{aligned}
	& \left( (\mu\times\mu) (S_{\varepsilon-2^{-t+1}}) - 3 \cdot 2^{-t} \right) \left( 1 - \delta \right) \\
	\leq & \underline{C}(G,x,\varepsilon,\omega,k) \\
	\leq & \overline{C}(G,x,\varepsilon,\omega,k) \\
	\leq & \left( (\mu\times\mu) (S_{\varepsilon+2^{-t+1}}) + 3 \cdot 2^{-t} \right) + \delta.
\end{aligned}
$$
Because $\delta$ is arbitrary, we conclude that given an $\omega \in \Sigma_m^+$, $k\in \mathbb{N}$ and any $\varepsilon \notin Q$, for any $t \in \mathbb{N}$, there exists a full measure subset $W_t \subseteq X$ such that any $x \in W_t$ satisfying
$$
\begin{aligned}
	& (\mu\times\mu) (S_{\varepsilon-2^{-t+1}}) - 3 \cdot 2^{-t} \\
	\leq & \underline{C}(G,x,\varepsilon,\omega,k) \\
	\leq & \overline{C}(G,x,\varepsilon,\omega,k) \\
	\leq & (\mu\times\mu) (S_{\varepsilon+2^{-t+1}}) + 3 \cdot 2^{-t}.
\end{aligned}
$$
We define $W(\omega, k, \varepsilon) := \bigcap_{t=1}^{+ \infty} W_t$. It is evident that $\mu (W(\omega, k, \varepsilon)) = 1$. Since $\varepsilon$ is a point of continuity of $\mu \times \mu(S_{\varepsilon})$, for any $x \in W(\omega, k, \varepsilon)$,
$$
(\mu\times\mu) (S_{\varepsilon}) = \underline{C}(G,x,\varepsilon,\omega,k) = \overline{C}(G,x,\varepsilon,\omega,k) = (\mu\times\mu) (S_{\varepsilon}),
$$
implying
$$
\lim_{n \rightarrow +\infty} C(G,x,\varepsilon,\omega,k,n) = (\mu\times\mu) (S_{\varepsilon})=\int_X \mu (B_{\omega,k}^G(x,\varepsilon))d\mu(x).
$$ \par
\end{proof}
\textbf{Step 3}. The proof proceeds by establishing the existence of a common full measure subset $W \subseteq X$ for any $\omega \in \Sigma_m^+$, $k\in \mathbb{N}$ and $\varepsilon \notin Q$.\par
\begin{proof}
In Step $2$, it is established that for any $\omega \in \Sigma_m^+, k \in \mathbb{N}, \varepsilon \notin Q$, there exists a full measure set $W(\omega, k, \varepsilon)$ such that for any $x \in W(\omega,k, \varepsilon)$,
$$
\lim_{n \rightarrow +\infty} C(G,x,\varepsilon,\omega,k,n) =\int_X \mu (B_{\omega,k}^G(x,\varepsilon))d\mu(x).
$$
Given that $\overline{h}_{cor}(G, \mu,2)$ and $\overline{h}_{cor}(G,x)$ are well-defined, we can substitute $\varepsilon \rightarrow 0$ with $\varepsilon_s \rightarrow 0$, where $\varepsilon_s \notin Q$ and $s \in \mathbb{N}$. This yields
$$
\begin{aligned}
\overline{h}_{cor}(G,\mu,2) & =\lim_{s \rightarrow +\infty} \limsup_{k \rightarrow +\infty} -\frac{1}{k} \int_{\Sigma_m^+} \log \left(\int_X \mu(B_{\omega,k}^G(x,\varepsilon_s)) d\mu(x)\right)d\mathbb{P}(\omega). \\
\overline{h}_{cor}(G,x)&=\lim_{s \rightarrow +\infty} \limsup_{k \rightarrow +\infty} - \frac{1}{k} \int_{\Sigma_m^+} \log \liminf_{n\rightarrow+\infty} C(G,x,\varepsilon_s,\omega,k,n) d\mathbb{P}(\omega).
\end{aligned}
$$
Moreover, $\Sigma_m^+$ is a compact metric space. For any $k,s \in \mathbb{N}$, consider the countable dense subset $\{\omega_r\}_{r=1}^{+\infty}$ of $\Sigma_m^+$, with $W(\omega_r,k,\varepsilon_s)$ representing the corresponding full measure subset of $X$. Define $W(k,s):=\bigcap_{r=1}^{+\infty}W(\omega_r,k,\varepsilon_s)$. For any $x \in W(k,s)$ and any $\omega \in \Sigma_m^+$, there exists $\omega_{r_0}$ such that $\omega_{r_0}|_{[1,k]}=\omega|_{[1,k]}$, where $\omega|_{[1,k]}$ denotes the initial $k$ elements of $\omega$. It is pertinent to note that $C(G,x,\varepsilon_s,\omega,k,n)$ is contingent upon $\omega|_{[1,k]}$, rather than $\omega$. For this reason, we have
$$
\begin{aligned}
\lim_{n \rightarrow +\infty} C(G,x,\varepsilon_s,\omega,k,n)&=\lim_{n \rightarrow +\infty} C(G,x,\varepsilon_s,\omega_{r_0},k,n) \\
&=\int_X \mu (B_{\omega_{r_0},k}^G(x,\varepsilon_s))d\mu(x)\\
&=\int_X \mu (B_{\omega,k}^G(x,\varepsilon_s))d\mu(x),
\end{aligned}
$$
implying
$$
\begin{aligned}
&- \frac{1}{k} \int_{\Sigma_m^+} \log \lim_{n \rightarrow +\infty} C(G,x,\varepsilon_s,\omega,k,n) d\mathbb{P}(\omega)\\
=&- \frac{1}{k}\int_{\Sigma_m^+} \log \left(\int_X \mu(B_{\omega,k}^G(x,\varepsilon_s)) d\mu(x)\right) d\mathbb{P}(\omega).
\end{aligned}
$$
Let $W:= \bigcap_{s,k=1}^{+\infty} W(k,s)$. It is obviously that $\mu(W) =1$. For any $x \in W$,
$$
\begin{aligned}
\overline{h}_{cor}(G,x)&=\lim_{s \rightarrow +\infty} \limsup_{k \rightarrow +\infty} - \frac{1}{k} \int_{\Sigma_m^+} \log \lim_{n\rightarrow+\infty} C(G,x,\varepsilon_s,\omega,k,n) d\mathbb{P}(\omega) \\
&=\lim_{s \rightarrow +\infty} \limsup_{k \rightarrow +\infty} -\frac{1}{k} \int_{\Sigma_m^+} \log \left(\int_X \mu(B_{\omega,k}^G(x,\varepsilon_s)) d\mu(x)\right)d\mathbb{P}(\omega) \\
&=\overline{h}_{cor}(G,\mu,2).
\end{aligned}
$$
Likewise, we can establish $\underline{h}_{cor}(G,\mu,2)=\underline{h}_{cor}(G,x)$ for $x \in W$. Hence, the Theorem \ref{th1.1} is demonstrated.
\end{proof}
\end{proof}

\begin{remark}
In \cite{PT}, the authors proved a similar result under stronger conditions, as stated below.
\begin{theorem}\cite{PT}
Let $(X,\mu)$ be a probability space, $G$ be a topological semigroup. If $\mu$ is $G$-ergodic and for any $\varphi \in L^1(\mu)$, there exists a full measure subset $Y \subseteq X$ such that for any $x \in Y$, the following holds
\begin{equation}\label{eq4.3}
\lim_{n \rightarrow + \infty}\frac{1}{n} \sum_{i=0}^{n-1} \varphi(f_{\omega,i}(x))=\int_{X}\varphi d\mu(x),
\end{equation}
and the convergence is uniform, then for any $x \in Y$, $\varepsilon >0$, $\omega \in \Sigma_m^+$ and $k \geq 1$, the equality
$$
\begin{aligned}
&\lim_{n \rightarrow + \infty}\frac{1}{n^2}\sharp \left\{ (i,j): 0\leq i,j \leq n-1, d_{\omega,k}^G(f_{\upsilon,i}(x),f_{\upsilon,j}(x))\leq \varepsilon \right\}\\
=&\int_{X}\mu(B_{\omega,k}^G(x,\varepsilon))d\mu(x)
\end{aligned}
$$
holds for any $\upsilon$.
\end{theorem}
In this paper, Theorem \ref{ergodic}(Theorem \ref{th1.1}) fails to meet the condition (\ref{eq4.3}) because full measure subset $\Omega_x \subseteq \Sigma_m^+$ depends on $x$ as noted in Remark \ref{re2.2}.
\end{remark}

\begin{corollary}
Let $G$ be a free semigroup acting on a compact metric space $X$, $(\Sigma_m^+, \mathbb{P})$ corresponding symbol space. For any $t \geq
1$, $(X,G^t)$ is the $t$-power system of $(X, G)$ and $(\Sigma_{m^t}^+, \mathbb{P}^t)$ is the corresponding symbol space. If $\mathbb{P}^t \times \mu$ is ergodic with respect to $F_t$, where $F_t$ is the skew product transformation acting on the $\Sigma_{m^t}^+ \times X$, defined as $F_t(\varpi,x) = (\sigma \varpi, g_{j_1}(x))$, $\varpi:= (j_1,j_2, \cdots) \in \Sigma_{m^t}^+$, then
$$
\overline{h}_{cor}(G^t,x)=t \cdot \overline{h}_{cor}(G,x), \quad
\underline{h}_{cor}(G^t,x)=t \cdot \underline{h}_{cor}(G,x),\qquad \mu-a.e. \quad x,
$$
\end{corollary}
\begin{proof}
Firstly, we assert that if $\mathbb{P}^t \times \mu$ is ergodic with respect to $F_t$, then $\mathbb{P}
\times \mu$ must also be ergodic with respect to $F$. To substantiate this claim, consider an invariant integrable function
$$
\varphi: \Sigma_m^+ \times X \longrightarrow \mathbb{R}.
$$
We can define a function $\psi: \Sigma_{m^t}^+ \times X \rightarrow \mathbb{R}$ as follows
$$
\psi(\varpi,x):= \varphi(\tau(\varpi),x),
$$
where $\tau : \Sigma_{m^t}^+ \rightarrow \Sigma_m^+$ is defined in Section 2.1. This function $\psi$ is integrable. For any $(\varpi,x)=((j_1,j_2,\cdots),x)$, the following transformation can be observed
$$
\begin{aligned}
\psi \circ F_t(\varpi,x)&=\psi((j_2,\cdots),g_{j_1}(x)) \\
&=\varphi(\sigma^t \tau(\varpi),f_{\tau(\varpi),t}(x)) \\
&=\varphi \circ F^t(\tau(\varpi),x) \\
&=\varphi(\tau(\varpi),x) \\
&=\psi(\varpi,x),
\end{aligned}
$$
meaning that $\psi$ is invariant. Since $\mathbb{P}^t \times \mu$ is ergodic, $\psi$ attains a constant value for $\mathbb{P}^t \times \mu$-almost everywhere $(\varpi, x)$. Furthermore, owing to the preservation of measure by both $\tau$ and $\tau^{-1}$, $\varphi$ takes on a constant value for $\mathbb{P} \times \mu$-almost everywhere $(\tau(\varpi),x)$. This observation implies the ergodicity of $\mathbb{P} \times \mu$ with respect to $F$. Utilizing Proposition 3.3
and Theorem 1.1, we derive the following equality
$$
\overline{h}_{cor}(G^t,x)=\overline{h}_{cor}(G^t,\mu,2)=t \cdot \overline{h}_{cor}(G,\mu,2)=t \cdot \overline{h}_{cor}(G,x), \quad \mu-a.e. \quad x.
$$
Similarly, the equality $\underline{h}_{cor}(G^t,x)=t \cdot \underline{h}_{cor}(G,x)$ for $\mu$-almost everywhere $x$ can be demonstrated using analogous reasoning.
\end{proof}

\leftline{\textbf{Problem}. In classical dynamical systems, \cite{S2018} established that}
$$
	\overline{h}_{cor}(f^t,x)=t \cdot \overline{h}_{cor}(f,x),  \quad
	\underline{h}_{cor}(f^t,x)=t \cdot \underline{h}_{cor}(f,x),
$$
holds for any $x$ via an innovative combinational approach. However, in the context of free semigroup actions, it remains unclear whether this power law persists for any $x$.

\section{Proof of Theorem \ref{th1.2}}
Prior to establishing Theorem \ref{th1.2} ($q=0$), we introduce a weak double-entropy condition for free semigroup actions, akin to the approach outlined in \cite{V}.
\begin{definition}
Let $G$ be a free semigroup acting on a compact metric space $X$, $\mu$ be a Borel probability measure on $X$. We say that $\mu$ satisfies the weak double-entropy condition of free semigroup actions if for sufficiently small $2\varepsilon$,, the
	$$
	\limsup_{k \rightarrow +\infty} \frac{1}{k} \log \sup_{x \in X} \frac{\mu(B_{\omega,k}^G(x,2\varepsilon))}{\mu(B_{\omega,k}^G(x,\varepsilon))}=0
	$$
	holds for almost everywhere $\omega \in \Sigma_m^+$.
\end{definition}

\begin{proof}[Proof of Theorem \ref{th1.2}]
Let $F(\omega, k,\varepsilon)$ be the $ (\omega, k,\varepsilon)$ spanning set with smallest cardinality. Since for any $x_i \in F(\omega, k, \varepsilon)$ and $x \in B^G_{\omega,k}(x_i,\varepsilon)$,
$$
B^G_{\omega,k}(x_i,\varepsilon) \subseteq B^G_{\omega,k}(x,2\varepsilon),
$$
we have
$$
\begin{aligned}
\int_X \mu \left(B^G_{\omega,k}(x,2\varepsilon) \right)^{-1} d \mu(x) \leq & \sum_{x_i \in F} \int_{B^G_{\omega ,k}(x_i,\varepsilon)} \mu \left(B^G_{\omega,k}(x,2\varepsilon) \right)^{-1} d\mu(x) \\
\leq & \sum_{x_i \in F} \mu \left(B^G_{\omega,k}(x_i,\varepsilon) \right)^{-1} \mu \left(B^G_{\omega,k}(x_i,\varepsilon) \right)\\
=&\sharp F(\omega,k,\varepsilon).
\end{aligned}
$$
Therefore
		$$
		\frac{1}{k} \int_{\Sigma_m^+} \log \left( \int_X \mu \left(B^G_{\omega,k}(x,2\varepsilon) \right)^{-1} d \mu(x) \right) d P(\omega) \leq \frac{1}{k} \int_{\Sigma_m^+} \log \sharp F(\omega,k,\varepsilon) d P(\omega).
		$$
		Taking $\limsup$ as $k \rightarrow +\infty$ and the limit as $\varepsilon \rightarrow 0$ on both sides of the inequality, we obtain
		$$
		\overline{h}_{cor}(G,\mu,0) \leq h_{top}(G).
		$$
Let $E(\omega,k,2\varepsilon)$ be the largest cardinality $ (\omega, k,2\varepsilon)$ separated set, with $2\varepsilon$ chosen sufficiently small such that $\mu$ satisfies the weak double-entropy condition of free semigroup actions. Consequently, there exist at most $\sharp E(\omega,k,2\varepsilon)$ pairwise disjoint Bowen ball $B^G_{\omega,k}(x_i,\varepsilon)$ on $X$, where $x_i \in E(\omega,k,2\varepsilon)$. For any $x \in B^G_{\omega,k}(x_i,\varepsilon)$, we have
$$
B^G_{\omega,k}(x,\varepsilon) \subseteq B^G_{\omega,k}(x_i,2\varepsilon).
$$
Thus,
		$$
		\begin{aligned}
			\int_X \mu \left(B^G_{\omega,k}(x,\varepsilon) \right)^{-1} d \mu(x) \geq& \sum_{x_i \in E} \int_{B^G_{\omega ,k}(x_i,\varepsilon)} \mu \left(B^G_{\omega,k}(x,\varepsilon) \right)^{-1} d\mu(x) \\
			\geq & \sum_{x_i \in E} \mu \left(B^G_{\omega,k}(x_i,2\varepsilon) \right)^{-1} \mu \left(B^G_{\omega,k}(x_i,\varepsilon) \right).
		\end{aligned}
		$$
Given that $\mu$ satisfies the weak double-entropy condition of free semigroup actions, we define $C_{\delta,K}$ for any $\delta > 0$ and $K \in \mathbb{N}$ as
		$$
		C_{\delta,K}:=\left\{\omega \in \Sigma_m^+: \frac{\mu(B_{\omega,k}^G(x,2\varepsilon))}{\mu(B_{\omega,k}^G(x,\varepsilon))} \leq e^{k\delta} \quad for \quad any \quad k>K \quad and \quad x \in X\right\}.
		$$
		It is clear that
		$$
		\begin{aligned}
		\left\{\omega \in \Sigma_m^+: \limsup_{k \rightarrow +\infty} \frac{1}{k} \log \sup_{x \in X} \frac{\mu(B_{\omega,k}^G(x,2\varepsilon))}{\mu(B_{\omega,k}^G(x,\varepsilon))}=0\right\} &=\bigcap_{\delta >0} \bigcup_{K=1}^{+ \infty} C_{\delta, K}\\
		&= \lim_{\delta \rightarrow 0}\lim_{K \rightarrow +\infty}C_{\delta,K}.
		\end{aligned}
		$$
		Therefore, for any $\eta >0$, there exists $\delta_0$ such that if $\delta < \delta_0$, then there exists $K_0=K_0(\delta)$ satisfying $\mathbb{P}(C_{\delta,K}) > 1- \eta$ holds for any $K > K_0$. For $\delta < \delta_0$ and $K > K_0$, consider $E(\omega,k,2\varepsilon)$ with $k > K$ and $2\varepsilon$ sufficiently small. For any $\omega \in C_{\delta,K}$, we have
		$$
		\begin{aligned}
		\int_X \mu \left(B^G_{\omega,k}(x,\varepsilon) \right)^{-1} d \mu(x) & \geq \sum_{x_i \in E} \mu \left(B^G_{\omega,k}(x_i,2\varepsilon) \right)^{-1} \mu \left(B^G_{\omega,k}(x_i,\varepsilon) \right)\\
		& \geq e^{-k\delta} \sharp E(\omega,k,2\varepsilon).
		\end{aligned}
		$$
		Consequently,
		$$
		\begin{aligned}
		&\frac{1}{k} \int_{\Sigma_m^+} \log \left( \int_X \mu \left(B^G_{\omega,k}(x,\varepsilon) \right)^{-1} d \mu(x) \right) d \mathbb{P}(\omega) \\
		 \geq & \frac{1}{k} \int_{C_{\delta,K}} \log \left( \int_X \mu \left(B^G_{\omega,k}(x,\varepsilon) \right)^{-1} d \mu(x) \right) d \mathbb{P}(\omega)\\
		 \geq & -\delta \mathbb{P}(C_{\delta,K}) + \frac{1}{k} \int_{C_{\delta,K}} \log \sharp E(\omega,k,2\varepsilon) d \mathbb{P}(\omega).
		\end{aligned}
		$$
		By taking $\liminf$ as $k \rightarrow +\infty$ and limit as $\varepsilon \rightarrow 0$ on both sides of the inequality, and employing Fatou's Lemma \cite{H}, we get
		$$
		\begin{aligned}
		& \lim_{\varepsilon \rightarrow 0} \liminf_{k \rightarrow +\infty} \frac{1}{k} \int_{\Sigma_m^+} \log \left( \int_X \mu \left(B^G_{\omega,k}(x,\varepsilon) \right)^{-1} d \mu(x) \right) d \mathbb{P}(\omega)\\
		 \geq & -\delta \mathbb{P}(C_{\delta,K}) + \int_{C_{\delta,K}} \lim_{\varepsilon \rightarrow 0} \liminf_{k \rightarrow +\infty} \frac{1}{k} \log \sharp E(\omega,k,2\varepsilon) d \mathbb{P}(\omega).
		\end{aligned}
		$$
		By the definition of topological entropy, we observe that for almost everywhere $\omega \in \Sigma_m^+$,
		$$
		\lim_{\varepsilon \rightarrow 0} \liminf_{k \rightarrow +\infty} \frac{1}{k} \log \sharp E(\omega,k,\varepsilon)=h_{top}(G).
		$$
		Therefore,
		$$
		\begin{aligned}
		&\lim_{\varepsilon \rightarrow 0} \liminf_{k \rightarrow +\infty} \frac{1}{k} \int_{\Sigma_m^+} \log \left( \int_X \mu \left(B^G_{\omega,k}(x,\varepsilon) \right)^{-1} d \mu(x) \right) d \mathbb{P}(\omega)\\
		\geq & -\delta \mathbb{P}(C_{\delta,K}) + h_{top}(G) \mathbb{P}(C_{\delta,K})\\
		\geq & -\delta + h_{top}(G) (1-\eta).
		\end{aligned}
		$$
		Since $\delta$ and $\eta$ are arbitrarily small, we conclude that
		$$
		\underline{h}_{cor}(G,\mu,0) \geq h_{top}(G).
		$$
\end{proof}

\begin{example}
\textcolor{red}{Let $\Sigma_2^+:=\{0,1\}^{\mathbb{N}}$ be a compact metric space, where the metric $d$ is defined as}
$$d((x_1,x_2,\cdots),(y_1,y_2,\cdots)):=2^{-\min\{i \geq 1: x_i \neq y_i\}},$$
and $\mu$ be a Bernoulli probability measure generated by $(\frac{1}{2},\frac{1}{2})$ on $(\Sigma_2^+,d)$. For any $x=(x_1,x_2,\cdots)$, the shift operator $f_1$ is defined as
$$
f_1(x):=(x_2,x_3,\cdots),
$$
and odometers $f_2$ (also known as adding machines) \cite{VO} is defined as
$$
f_2(x):=x+(1,0,0,\cdots)=(y_1,y_2,y_3,\cdots),
$$
where $(y_1,y_2,y_3,\cdots)$ is determined by the following process.\\
If $x_1+1=1$, then $y_1=1$ and $\delta_2=0$,\\
if $x_1+1=2$, then $y_1=0$ and $\delta_2=1$.\\
For every $n \geq 2$,\\
if $x_n+\delta_n=1$, then $y_n=1$ and $\delta_{n+1}=0$,\\
if $x_n+\delta_n=2$, then $y_n=0$ and $\delta_{n+1}=1$.\\
Denote $i=min\{j \geq 1: x_j=0\}$, that is, $x=(1, 1,\cdots, 1, 0,x_{i+1},x_{i+2},\cdots)$. Hence, we get a simple expression of $f_2$ as follows
$$
f_2(x):=x+(1,0,0,\cdots)=(0,\cdots,0,1,x_{i+1},x_{i+2},\cdots).
$$
It is noted that if $i=+\infty$, that is, $x=(1,1,\cdots)$, then $f_2(x)=(0,0,\cdots)$. $G$ is the free semigroup generated by $\{f_1,f_2\}$ acting on compact metric space $(\Sigma_2^+,d)$. $(\Sigma_2^+,\mathbb{P})$ is the corresponding symbol space where $\mathbb{P}$ is the Bernoulli probability measure generated by $(\frac{1}{2},\frac{1}{2})$. For the sake of convenience, let $\varepsilon = 2^{-t}$, $t \in \mathbb{N}$. Note that
$$
f_1^{-1}B(f_1(x),\varepsilon)=_{2}[x_2,x_3,\cdots,x_{t+1}]_{t+1}, \quad \quad f_2^{-1}B(f_2(x),\varepsilon)=_{1}[x_1,x_2,\cdots,x_t]_{t},
$$
where cylinder
$$
_{a}[i_1,i_2,\cdots,i_j]_{a+j-1}:=\{y=(y_1,y_2,\cdots) \in \Sigma_2^+: y_{a+t}=i_{t+1}, 0 \leq t \leq j-1 \}.
$$
We claim that for any $\omega:=(i_1,i_2,\cdots) \in \Sigma_2^+$, $k \geq 1$ and $\varepsilon =2^{-t}$, it is verified that
$$
B_{\omega,k}^G(x,\varepsilon)=_{1}[x_1,x_2,\cdots,x_{t+s_{\omega,k}}]_{t+s_{\omega,k}}
$$
where $s_{\omega,k}:=\sharp \{1 \leq j \leq k-1: i_j =1 \}$. Next we provide the proof of the claim.\par
It is known that
$$
B_{\omega,k}^G(x,\varepsilon):=B(x,\varepsilon) \bigcap \cdots \bigcap f_{\omega,k-2}^{-1}B(f_{\omega,k-2}(x),\varepsilon)\bigcap f_{\omega,k-1}^{-1}B(f_{\omega,k-1}(x),\varepsilon).
$$
If $f_{i_{k-1}}=f_1$, then $f_{\omega,k-1}^{-1}B(f_{\omega,k-1}(x),\varepsilon)=f_{\omega,k-2}^{-1} \circ f_1^{-1} B(f_1 \circ f_{\omega,k-2}(x),\varepsilon)$, implying that
$$
\begin{aligned}
B_{\omega,k}^G(x,\varepsilon):=&B(x,\varepsilon) \bigcap f_{\omega,1}^{-1}B(f_{\omega,1}(x),\varepsilon) \bigcap \cdots \bigcap f_{\omega,k-2}^{-1}B(f_{\omega,k-2}(x),\varepsilon)\\ &\bigcap f_{\omega,k-2}^{-1} \circ f_1^{-1} B(f_1 \circ f_{\omega,k-2}(x),\varepsilon)\\
=&B(x,\varepsilon) \bigcap f_{\omega,1}^{-1}B(f_{\omega,1}(x),\varepsilon) \bigcap \cdots \bigcap f_{\omega,k-3}^{-1}B(f_{\omega,k-3}(x),\varepsilon) \\
&\bigcap f_{\omega,k-2}^{-1}\left(B(f_{\omega,k-2}(x),\varepsilon) \bigcap f_1^{-1} B(f_1 \circ f_{\omega,k-2}(x),\varepsilon)\right) \\
=&B(x,\varepsilon) \bigcap f_{\omega,1}^{-1}B(f_{\omega,1}(x),\varepsilon) \bigcap \cdots \bigcap f_{\omega,k-3}^{-1}B(f_{\omega,k-3}(x),\varepsilon) \\
&\bigcap f_{\omega,k-2}^{-1}B(f_{\omega,k-2}(x),\frac{\varepsilon}{2}).
\end{aligned}
$$
If $f_{i_{k-1}}=f_2$, then $f_{\omega,k-1}^{-1}B(f_{\omega,k-1}(x),\varepsilon)=f_{\omega,k-2}^{-1} \circ f_2^{-1} B(f_2 \circ f_{\omega,k-2}(x),\varepsilon)$, implying that
$$
B_{\omega,k}^G(x,\varepsilon):=B(x,\varepsilon) \bigcap f_{\omega,1}^{-1}B(f_{\omega,1}(x),\varepsilon) \bigcap \cdots \bigcap f_{\omega,k-2}^{-1}B(f_{\omega,k-2}(x),\varepsilon).
$$
Based on these consideration, we establish the claim through induction. \par
For the base case $k=1$, we have $B_{\omega,k}^G(x,\varepsilon)=B(x,\varepsilon)=_{1}[x_1,x_2,\cdots,x_t]_{t}$. \par
Now, for $k=2$, if $f_{i_1}=f_1$, then $B_{\omega,k}^G(x,\varepsilon)=B(x,\frac{\varepsilon}{2})=_{1}[x_1,x_2,\cdots,x_{t+1}]_{t+1}$ and if $f_{i_1}=f_2$, then $B_{\omega,k}^G(x,\varepsilon)=B(x,\varepsilon)=_{1}[x_1,x_2,\cdots,x_t]_{t}$.\par
Now, Assuming that the assertion holds for $k-1$, we examine the case $k$. Let us define
$$
s_{\omega,k-1}:=\sharp \{1 \leq j \leq k-2: i_j =1 \},  \qquad s_{\omega,k}:=\sharp \{1 \leq j \leq k-1: i_j =1 \}.
$$
If $f_{i_{k-1}}=f_1$, meaning $s_{\omega,k}=s_{\omega,k-1}+1$, we have
$$
\begin{aligned}
B_{\omega,k}^G(x,\varepsilon)=&B(x,\varepsilon) \bigcap f_{\omega,1}^{-1}B(f_{\omega,1}(x),\varepsilon) \bigcap \cdots \bigcap f_{\omega,k-3}^{-1}B(f_{\omega,k-3}(x),\varepsilon) \\
&\bigcap f_{\omega,k-2}^{-1}B(f_{\omega,k-2}(x),\frac{\varepsilon}{2}).
\end{aligned}
$$
When $f_{i_{k-2}}=f_1$,
$$
\begin{aligned}
&f_{\omega,k-3}^{-1}B(f_{\omega,k-3}(x),\varepsilon) \bigcap f_{\omega,k-2}^{-1}B(f_{\omega,k-2}(x),\frac{\varepsilon}{2})\\
=&f_{\omega,k-3}^{-1}\left( B(f_{\omega,k-3}(x),\varepsilon) \bigcap f_1^{-1}B(f_1 \circ f_{\omega,k-3}(x),\frac{\varepsilon}{2}) \right)\\
=&f_{\omega,k-3}^{-1} B(f_{\omega,k-3}(x),\frac{\varepsilon}{2^2})\\
\subseteq &f_{\omega,k-3}^{-1}B(f_{\omega,k-3}(x),\frac{\varepsilon}{2}) \bigcap f_{\omega,k-2}^{-1}B(f_{\omega,k-2}(x),\frac{\varepsilon}{2})\\
\subseteq &f_{\omega,k-3}^{-1}B(f_{\omega,k-3}(x),\varepsilon) \bigcap f_{\omega,k-2}^{-1}B(f_{\omega,k-2}(x),\frac{\varepsilon}{2}).
\end{aligned}
$$
When $f_{i_{k-2}}=f_2$,
$$
\begin{aligned}
&f_{\omega,k-3}^{-1}B(f_{\omega,k-3}(x),\varepsilon) \bigcap f_{\omega,k-2}^{-1}B(f_{\omega,k-2}(x),\frac{\varepsilon}{2})\\
=&f_{\omega,k-3}^{-1}B(f_{\omega,k-3}(x),\varepsilon) \bigcap f_{\omega,k-3}^{-1}B(f_{\omega,k-3}(x),\frac{\varepsilon}{2})\\
\subseteq &f_{\omega,k-3}^{-1}B(f_{\omega,k-3}(x),\frac{\varepsilon}{2}) \bigcap f_{\omega,k-2}^{-1}B(f_{\omega,k-2}(x),\frac{\varepsilon}{2})\\
\subseteq &f_{\omega,k-3}^{-1}B(f_{\omega,k-3}(x),\varepsilon) \bigcap f_{\omega,k-2}^{-1}B(f_{\omega,k-2}(x),\frac{\varepsilon}{2}).
\end{aligned}
$$
Consequently,
$$
\begin{aligned}
B_{\omega,k}^G(x,\varepsilon)=&B(x,\varepsilon) \bigcap f_{\omega,1}^{-1}B(f_{\omega,1}(x),\varepsilon) \bigcap \cdots \bigcap f_{\omega,k-3}^{-1}B(f_{\omega,k-3}(x),\frac{\varepsilon}{2}) \\
&\bigcap f_{\omega,k-2}^{-1}B(f_{\omega,k-2}(x),\frac{\varepsilon}{2}).
\end{aligned}
$$
Iterating this process $k-2$ times yields
$$
\begin{aligned}
B_{\omega,k}^G(x,\varepsilon)&=B_{\omega,k-1}^G(x,\frac{\varepsilon}{2})\\
&=_{1}[x_1,x_2,\cdots,x_{t+s_{\omega,k-1}+1}]_{t+s_{\omega,k-1}+1}\\
&=_{1}[x_1,x_2,\cdots,x_{t+s_{\omega,k}}]_{t+s_{\omega,k}}.
\end{aligned}
$$
If $f_{i_{k-1}}=f_2$, that is $s_{\omega,k}=s_{\omega,k-1}$, we obtain
$$
\begin{aligned}
B_{\omega,k}^G(x,\varepsilon):=&B(x,\varepsilon) \bigcap f_{\omega,1}^{-1}B(f_{\omega,1}(x),\varepsilon) \bigcap \cdots \bigcap f_{\omega,k-2}^{-1}B(f_{\omega,k-2}(x),\varepsilon)\\
=&B_{\omega,k-1}^G(x,\varepsilon)\\
=&_{1}[x_1,x_2,\cdots,x_{t+s_{\omega,k-1}}]_{t+s_{\omega,k-1}}\\
=&_{1}[x_1,x_2,\cdots,x_{t+s_{\omega,k}}]_{t+s_{\omega,k}}.
\end{aligned}
$$
Hence, the claim is demonstrated. Based on this claim, $\mu$ satisfies weak double-entropy condition of free semigroup actions. By Theorem\ref{th1.2},
$$
\begin{aligned}
h_{top}(G)=& \lim_{\varepsilon \rightarrow 0}\lim_{k \rightarrow +\infty}\frac{1}{k}\int_{\Sigma_2^+}\log\left(\int_{X}\mu(B_{\omega,k}^G(x,\varepsilon))^{-1}d\mu(x)\right)d\mathbb{P}(\omega)\\
=&\lim_{\varepsilon \rightarrow 0}\lim_{k \rightarrow +\infty}\frac{1}{k}\int_{\Sigma_2^+}\log\left(\int_{X}\mu(_{1}[x_1,x_2,\cdots,x_{t+s_{\omega,k}}]_{t+s_{\omega,k}})^{-1}d\mu(x)\right)d\mathbb{P}(\omega)\\
=&\lim_{\varepsilon \rightarrow 0}\lim_{k \rightarrow +\infty}\frac{1}{k}\int_{\Sigma_2^+}\log\left(\int_{X} 2^{t+s_{\omega,k}}d\mu(x)\right)d\mathbb{P}(\omega)\\
=&\lim_{\varepsilon \rightarrow 0}\lim_{k \rightarrow +\infty}\frac{1}{k}\int_{\Sigma_2^+}\log 2^{t+s_{\omega,k}}d\mathbb{P}(\omega)\\
=&\lim_{\varepsilon \rightarrow 0}\lim_{k \rightarrow +\infty}\frac{1}{k}\sum_{s=0}^{k-1}C_{k-1}^s 2^{-(k-1)}\log2^{t+s}\\
=&\lim_{\varepsilon \rightarrow 0}\lim_{k \rightarrow +\infty} \frac{\log2}{2^{k-1}k}\sum_{s=0}^{k-1}C_{k-1}^s (t+s)\\
=&\lim_{\varepsilon \rightarrow 0}\lim_{k \rightarrow +\infty} \frac{\log2}{2^{k-1}k} \left(\sum_{s=0}^{k-1}C_{k-1}^s t+\sum_{s=0}^{k-1}C_{k-1}^s s\right)\\
=&\lim_{\varepsilon \rightarrow 0}\lim_{k \rightarrow +\infty} \frac{\log2}{2^{k-1}k} \left(t2^{k-1}+(k-1)2^{k-2}\right)\\
=&\lim_{\varepsilon \rightarrow 0}\lim_{k \rightarrow +\infty} \left(\frac{t\log2}{k}+\frac{\log2}{2}\frac{k-1}{k}\right)\\
=&\frac{\log2}{2}.
\end{aligned}
$$
\end{example}

\section{Proofs of Theorem \ref{th1.3}, \ref{th1.4}, \ref{th1.5}}
Motivated by Theorem \ref{th2.3}, we introduce the notions of the lower(upper) local entropy of free semigroup actions as follows.
\begin{definition}
Let $G$ be a free semigroup acting on a compact metric space $X$, $\mu$ a Borel probability measure on $X$. The lower(upper) local entropy of free semigroup actions is defined as
	$$
	\begin{aligned}
		& h(\omega,x):= \lim_{\varepsilon \rightarrow 0} \liminf_{k \rightarrow +\infty} -\frac{1}{k} \log \mu (B_{\omega,k}^G(x,\varepsilon)), \\
		& H(\omega,x):= \lim_{\varepsilon \rightarrow 0} \limsup_{k \rightarrow +\infty} -\frac{1}{k} \log \mu (B_{\omega,k}^G(x,\varepsilon)).
	\end{aligned}
    $$
\end{definition}
We define the limit process of $h(\omega,x)$ to be uniformly with respect to $x$ for almost everywhere $\omega$ if, for almost everywhere $\omega$, there exists a full measure subset $A(\omega) \subseteq X$ such that for any $\delta > 0$, there exists $\varepsilon_0:=\varepsilon_0(\omega,\delta)$ such that for $\varepsilon \leq \varepsilon_0$, there exists $K:=K(\omega,\delta,\varepsilon)$ such that if $k > K$, then
$$
h(\omega, x) - \delta \leq -\frac{1}{k} \log \mu(B^G_{\omega, k}(x, \varepsilon))
$$
holds for any $x \in A(\omega)$. We define he limit process of $h(\omega,x)$ to be uniformly with respect to $(\omega,x)$ if there exists a full measure subset $A \subseteq \Sigma_m^+ \times X$ such that for any $\delta > 0$, there exists $\varepsilon_0:=\varepsilon_0(\delta)$ such that for $\varepsilon \leq \varepsilon_0$, there exists $K:=K(\delta,\varepsilon)$ such that if $k > K$, then
$$
h(\omega, x) - \delta \leq -\frac{1}{k} \log \mu(B^G_{\omega, k}(x, \varepsilon))
$$
holds for any $(\omega,x) \in A$.\par
Prior to establishing Theorem \ref{th1.3} ($q \geq 1$), we require the following lemma.
\begin{lemma}\label{lemma6.1}
Let $G$ be a free semigroup acting on a compact metric space $X$, $\mu$ a $G$-invariant Borel probability measure on $X$. Then
$$
\overline{h}_{cor}(G,\mu,1) \leq h_{\mu}(G).
$$
\end{lemma}
\begin{proof}
For any $\varepsilon >0$, there exists a finite partition $\xi$ such that ${\rm diam}(\xi) \leq \varepsilon$. Considering that
$D(\bigvee_{i=0}^{k-1} f_{\omega,i}^{-1} \xi, x) \subseteq B_{\omega, k}^G(x,\varepsilon)$ where $D(\bigvee_{i=0}^{k-1} f_{\omega,i}^{-1} \xi, x)$ represents the element of $\bigvee_{i=0}^{k-1} f_{\omega,i}^{-1} \xi$ containing $x$, we obtain
\begin{equation}\label{eq6.1}
\begin{aligned}
&-\frac{1}{k} \int_{\Sigma_m^+} \int_X \log \mu (D(\bigvee_{i=0}^{k-1} f_{\omega,i}^{-1} \xi, x)) d\mu(x) d\mathbb{P}(\omega)\\
\geq &-\frac{1}{k} \int_{\Sigma_m^+} \int_X \log \mu \left( B_{\omega, k}^G(x,\varepsilon) \right) d\mu(x)d\mathbb{P}(\omega).
\end{aligned}
\end{equation}
By definition of $H_{\mu}(\bigvee_{i=0}^{k-1} f_{\omega, i}^{-1} \xi)$, we get
$$
\frac{1}{k}\int_{\Sigma_m^+} H_{\mu}(\bigvee_{i=0}^{k-1} f_{\omega, i}^{-1} \xi) d\mathbb{P}(\omega) \geq -\frac{1}{k}
\int_{\Sigma_m^+} \int_X \log \mu\left( B_{\omega, k}^G(x,\varepsilon)\right) d\mu(x) d\mathbb{P}(\omega).
$$
Taking $\limsup$ as $k \rightarrow +\infty$ and limit as $\varepsilon \rightarrow 0$, we have
$$
h_{\mu}(G) \geq h_{\mu}(G, \xi) \geq \overline{h}_{cor}(G,\mu,1).
$$
\end{proof}
\begin{remark}
	Lemma \ref{lemma6.1} generalizes the result of E.Verbitskiy \cite{V} to the free semigroup actions.
\end{remark}
\begin{proof}[Proof of Theorem \ref{th1.3}]
According to Theorem 2.3., there exists measurable set $A \subseteq \Sigma_m^+ \times X$ with full measure, such that for any $(\omega,x) \in
A$, $h(\omega,x)=h_{\mu}(G)$. Denote
$$
\mathcal{P}(A):=\{\omega: \exists x, s.t. (\omega,x) \in A\}, \qquad A(\omega):=\{x: (\omega,x) \in A\}.
$$
Since $\mathbb{P} \times \mu(A)=1$, for almost everywhere $\omega \in \mathcal{P}(A)$, $\mu(A(\omega))=1$.
Under the given assumption, we can choose $A(\omega)$ to satisfy the following condition. For almost everywhere $\omega \in \mathcal{P}(A)$ and any $\delta > 0$, there exists $\varepsilon_0:=\varepsilon_0(\omega,\delta)$, such that for any $\varepsilon \leq \varepsilon_0$, there exists $K:=K(\omega,\delta,\varepsilon)$ such that if $k > K$, then
$$
h(\omega, x) - \delta \leq -\frac{1}{k} \log \mu(B^G_{\omega, k}(x, \varepsilon))
$$
holds for any $x \in A(\omega)$.
Similar to Theorem \ref{th1.2}, for any $\delta > 0$, $\varepsilon > 0$ and $K \in \mathbb{N}$, we define
$$
\begin{aligned}
\Omega_{\delta,\varepsilon,K}:=&\Big\{\omega \in \mathcal{P}(A): h_{\mu}(G) - \delta \leq -\frac{1}{k} \log \mu(B^G_{\omega, k}(x, \varepsilon)) \ for \ any \ k >K,  \ x \in A(\omega)\Big\}.
\end{aligned}
$$
As
$$
\lim_{\delta \rightarrow 0}\lim_{\varepsilon \rightarrow 0} \lim_{K \rightarrow +\infty} \mathbb{P}(\Omega_{\delta,\varepsilon,K})=1,
$$
for any $\eta >0$, $\delta > 0$. there exist $\varepsilon$, $K$ such that $\mathbb{P}(\Omega_{\delta,\varepsilon,K}) > 1-\eta$. Given $k > K$ and $q >1$, we have
$$
\begin{aligned}
&-\frac{1}{k} \frac{1}{q-1} \int_{\Sigma_m^+} \log \left(\int_X \mu(B_{\omega,k}^G(x,\varepsilon))^{q-1} d\mu(x)\right) d\mathbb{P}(\omega) \\
\geq &-\frac{1}{k} \frac{1}{q-1} \int_{\Omega_{\delta,\varepsilon,K}} \log \left(\int_{A(\omega)} e^{-k(q-1)(h_{\mu}(G)-\delta)} d\mu(x)\right)
d\mathbb{P}(\omega) \\
=& (h_{\mu}(G)-\delta) \mathbb{P}(\Omega_{\delta,\varepsilon,K}).\\
\end{aligned}
$$
If $h_{\mu}(G)=0$, then
$$
-\frac{1}{k} \frac{1}{q-1} \int_{\Sigma_m^+} \log \left(\int_X \mu(B_{\omega,k}^G(x,\varepsilon))^{q-1} d\mu(x)\right) d\mathbb{P}(\omega) \geq -\delta.
$$
If $h_{\mu}(G)>0$, we can assume $h_{\mu}(G)-\delta > 0$ owing to $\delta$ is arbitrary, then
$$
-\frac{1}{k} \frac{1}{q-1} \int_{\Sigma_m^+} \log \left(\int_X \mu(B_{\omega,k}^G(x,\varepsilon))^{q-1} d\mu(x)\right) d\mathbb{P}(\omega) \geq (h_{\mu}(G)-\delta)(1-\eta).
$$
As $\delta$ and $\eta$ are arbitrary, taking $\liminf$ as $k\rightarrow +\infty$ and limit as $\varepsilon \rightarrow 0$ on both sides yields
$$
\underline{h}_{cor}(G,\mu,q) \geq h_{\mu}(G).
$$
By utilizing Lemma 6.1 and Proposition 3.2 (4), we can establish for any $q \geq 1$ that
$$
h_{\mu}(G)=\overline{h}_{cor}(G,\mu,q)=\underline{h}_{cor}(G,\mu,q).
$$
In particularly, based on Theorem 1.1,
$$
h_{\mu}(G)=h_{cor}(G,\mu,2)=h_{cor}(G,x), \qquad \mu-a.e. \quad x \in X.
$$
\end{proof}
\begin{remark}
Following the approach in \cite{V}, we establish $h_{\mu}(G) = h_{cor}(G,\mu,1)$ without the requirement concerning $h(\omega,x)$, as stated in Proposition \ref{prop6.1}.
\end{remark}

\begin{proposition}\label{prop6.1}
\textcolor{red}{Let $G$ be a free semigroup acting on a compact metric space $X$ and $\mu$ be a Borel probability measure on $X$ such that $\mathbb{P} \times \mu$ is ergodic with respect to skew product transformation $F$. If $\mu$ satisfies $h_{\mu}(G) < +\infty$, then}
$$
h_{cor}(G,\mu,1)=h_{\mu}(G).
$$
\end{proposition}
	\begin{proof}
Given $h_{\mu}(G) < +\infty$ and the ergodicity of $\mathbb{P} \times \mu$, it follows that $h(\omega,x)=h_{\mu}(G)$ for almost everywhere $(\omega,x)$. By Fatou's Lemma \cite{H}, we obtain
$$
\begin{aligned}
&\lim_{\varepsilon \rightarrow 0} \liminf_{k \rightarrow +\infty} -\frac{1}{k} \int_{\Sigma_m^+} \int_X \log \mu \left( B_{\omega, k}^G(x,\varepsilon) \right) d\mu(x)d\mathbb{P}(\omega)\\
\geq & \int_{\Sigma_m^+} \int_X \lim_{\varepsilon \rightarrow 0} \liminf_{k \rightarrow +\infty} -\frac{1}{k} \log \mu \left( B_{\omega, k}^G(x,\varepsilon) \right) d\mu(x)d\mathbb{P}(\omega)\\
= & h_{\mu}(G).
\end{aligned}
$$
Thus, $\underline{h}_{cor}(G,\mu,1) \geq h_{\mu}(G)$. By utilizing Lemma \ref{lemma6.1}, we conclude that
$$
h_{cor}(G,\mu,1) = h_{\mu}(G).
$$
\end{proof}
We now turn our attention to the scenario where $0 \leq q \leq 1$, as outlined in Theorem \ref{th1.4}.
\begin{proof}[Proof of Theorem \ref{th1.4}]
There exists a set $ A \subseteq \Sigma_m^+ \times X$ of full measure such that for any $ (\omega,x) \in A$, $h(\omega,x)$ exists and $h(\omega,x) \geq h_{top}(G)$ by assumption. We define
$$
\mathcal{P}(A):=\{\omega: \exists x, s.t. (\omega,x) \in A\}, \qquad A(\omega):=\{x: (\omega,x) \in A\}.
$$
By Fatou's lemma \cite{H}, we obtain
	$$
	\begin{aligned}
		&\lim_{\varepsilon \rightarrow 0} \liminf_{k \rightarrow +\infty} -\frac{1}{k} \int_{\Sigma_m^+} \int_X \log \mu (B_{\omega,k}^G(x,\varepsilon)) d\mu (x) d \mathbb{P} (\omega) \\
		\geq & \int_{\Sigma_m^+} \int_X \lim_{\varepsilon \rightarrow 0} \liminf_{k \rightarrow +\infty} -\frac{1}{k} \log \mu (B_{\omega,k}^G(x,\varepsilon)) d\mu (x) d \mathbb{P} (\omega) \\
		= & \int_{\mathcal{P}(A)} \int_{A(\omega)} \lim_{\varepsilon \rightarrow 0} \liminf_{k \rightarrow +\infty} -\frac{1}{k} \log \mu (B_{\omega,k}^G(x,\varepsilon)) d\mu (x) d \mathbb{P} (\omega) \\
		\geq & h_{top}(G).
	\end{aligned}
	$$
Thus, we have
$$
\underline{h}_{cor}(G,\mu,1) \geq h_{top}(G).
$$
On the other hand, from the proof of Theorem \ref{th1.2}, it follows that
$$
\overline{h}_{cor}(G,\mu,0) \leq h_{top}(G).
$$
Utilizing Proposition 3.2 (4), we conclude that for any $0 \leq q \leq 1$,
$$
h_{top}(G)=h_{cor}(G,\mu,q).
$$
Referencing \cite{B1}, it follows that
\begin{equation}\label{eq6.2}
h_{top}(G)= \sup \{h_{\mu}(G): \mu \in M(X,G)\},
\end{equation}
where $M(X,G)$ denotes the space of $G$-invariant probability measure of $(X,G)$. By employing Lemma 6.1, we establish
$$
h_{top}(G)=h_{cor}(G,\mu,1)= h_{\mu}(G),
$$
indicating that $\mu$ represents the measure of maximum entropy.
\end{proof}
Finally, we demonstrate Theorem \ref{th1.5} for $q \leq 0$.
\begin{proof}[Proof of Theorem \ref{th1.5}]
Under the given assumption, there exists a full measure subset $A \subseteq \Sigma_m^+ \times X$ such that for any $\delta > 0$, there exists $\varepsilon_0:=\varepsilon_0(\delta)$, such that for any $\varepsilon \leq \varepsilon_0$, there exists $K:=K(\delta,\varepsilon)$ satisfying if $k > K$, then
$$
\mu(B^G_{\omega, k}(x, \varepsilon))^{q-1} \leq e^{-k(q-1)(h_{\mu}(G)+\delta)}
$$
holds for any $(\omega,x) \in A$ and $q \leq 0$. Consequently, for any $q \leq 0$, $\varepsilon \leq \varepsilon_0$ and $k > K$, we derive
$$
\begin{aligned}
&-\frac{1}{k} \frac{1}{q-1} \int_{\Sigma_m^+} \log \left(\int_X \mu(B_{\omega,k}^G(x,\varepsilon))^{q-1} d\mu(x)\right) d\mathbb{P}(\omega) \\
\leq &-\frac{1}{k} \frac{1}{q-1} \int_{\mathcal{P}(A)} \log \left(\int_{A(\omega)} e^{-k(q-1)(h_{\mu}(G)+\delta)} d\mu(x)\right)
d\mathbb{P}(\omega) \\
=& (h_{\mu}(G)+\delta).\\
\end{aligned}
$$
Taking the $\limsup$ as $k \rightarrow +\infty$ and the limit as $\varepsilon \rightarrow 0$ on both sides, we obtain for any $q \leq 0$,
$$
\overline{h}_{cor}(G,\mu,q) \leq h_{\mu}(G).
$$
By employing Theorem \ref{th1.2} and equality (\ref{eq6.2}), we conclude that for any $q \leq 0$,
$$
h_{cor}(G,\mu,q)=h_{\mu}(G).
$$
Specifically,
$$
h_{top}(G) =h_{cor}(G,\mu,0)= h_{\mu}(G),
$$
signifying that $\mu$ represents the measure of maximum entropy.
\end{proof}

\begin{example}
Let $X$ be a compact metric group with the right-invariant metric $d$ and Haar measure $\mu$, which also exhibits right-invariance property (\cite{MZ}). Consider $L: X \rightarrow X$ as a group automorphism. We define a continuous map $f_i: X \rightarrow X$ for each $x_i$ in a finite set $\{x_i: 1 \leq i \leq m\} \subseteq X$ as $f_i(x):=L(x)\cdot x_i$. Given automorphism nature of $L$, the inverse map of $f_i$ is denoted as $f_i^{-1}(x):=L^{-1}(x)\cdot L^{-1}(x_i^{-1})$. Let $G$ be the free semigroup generated by $\{f_i: 1 \leq i \leq m\}$. Notably, due to the right-invariance of metric $d$, for any $\omega=(i_1,i_2,\cdots)$ and $j\geq 0$, we have
$$
B(f_{\omega,j}(x),\varepsilon)=B(e,\varepsilon)\cdot f_{\omega,j}(x)
$$
where $e$ represents the identity element. Consequently,
$$
f_{\omega,j}^{-1}B(f_{\omega,j}(x),\varepsilon)=f_{\omega,j}^{-1}\left(B(e,\varepsilon)\cdot f_{\omega,j}(x) \right).
$$
Notice that for any $z \in B(e,\varepsilon)$,
$$
\begin{aligned}
&f_{\omega,j}^{-1}\left(z\cdot f_{\omega,j}(x) \right)\\
=&f_{i_1}^{-1}\circ f_{i_2}^{-1} \circ \cdots \circ f_{i_j}^{-1} \left( z\cdot f_{\omega,j}(x)\right) \\
=&f_{i_1}^{-1}\circ f_{i_2}^{-1} \circ \cdots \circ f_{i_{j-1}}^{-1} L^{-1}\left( z\cdot f_{\omega,j}(x)\right) \cdot L^{-1}(x_{i_j}^{-1})\\
=&f_{i_1}^{-1}\circ f_{i_2}^{-1} \circ \cdots \circ f_{i_{j-1}}^{-1} L^{-1}(z)L^{j-1}(x)L^{j-2}(x_{i_1})\cdots x_{i_{j-1}}\cdot L^{-1}(x_{i_j})\cdot L^{-1}(x_{i_j}^{-1})\\
=&f_{i_1}^{-1}\circ f_{i_2}^{-1} \circ \cdots \circ f_{i_{j-1}}^{-1} L^{-1}(z)L^{j-1}(x)L^{j-2}(x_{i_1})\cdots x_{i_{j-1}}\\
&\cdots \\
=&L^{-j}(z) \cdot x.
\end{aligned}
$$
Hence,
$$
f_{\omega,j}^{-1}B(f_{\omega,j}(x),\varepsilon)=f_{\omega,j}^{-1}\left(B(e,\varepsilon)\cdot f_{\omega,j}(x) \right)=L^{-j}\left(B(e,\varepsilon)\right)\cdot x,
$$
which implies that that $B_{\omega,k}^G(x,\varepsilon)=\left(\bigcap_{j=0}^{k-1}L^{-j}\left( B(e,\varepsilon)\right)\right)\cdot x$. Moreover, since $\mu$ is right-invariant, it follows that the limit processes of $h(\omega,x)$ and $H(\omega,x)$ are uniform. In classical dynamical systems, Bowen \cite{B1971} introduced the concept of homogeneous measure, where the Haar measure serves as a homogeneous measure when $L: X \rightarrow X$ denotes an affine transformation acting on a compact metric group $X$. Verbitskiy \cite{V} computed the order correlation entropy of homogeneous measure and demonstrated the equality between any order correlation entropy, measure-theoretic entropy, and topological entropy. This equivalence extends to free semigroup actions, where these entropies, computed by definition, also coincide due to the inherent strong uniformity of homogeneous measures.
Theorems \ref{th1.3}, \ref{th1.4}, and \ref{th1.5} are established by relaxing the uniformity constraints on the measure while incorporating additional conditions.
\end{example}

\section*{Declarations}
\leftline{\textbf{Conflict of interest} The authors declare their is no conflict of interest.}

\section*{Acknowledgments}
The authors express their sincere appreciation for the insightful remarks and constructive suggestions provided by the referees, which have significantly enhanced the quality of this manuscript. Additionally, the authors would like to acknowledge Prof. Xiaogang Lin for his guidance in English academic writing.

\section*{References}
\bibliographystyle{plain}
\bibliography{references}

\end{document}